\documentclass{paper}
\RequirePackage{graphicx}% images

\usepackage{here}        % wrap figures/tables in text (i.e., Di Vinci style)
\usepackage{wrapfig}        % wrap figures/tables in text (i.e., Di Vinci style)
\usepackage{subfigure}      % subcaptions for subfigures
\usepackage{subfigmat}      % matrices of similar subfigures, aka small multiples

\usepackage{varioref}       % smart page, figure, table, and equation referencing
\usepackage{wrapfig}        % wrap figures/tables in text (i.e., Di Vinci style)
\usepackage{threeparttable} % tables with footnotes
\usepackage{dcolumn}        % decimal-aligned tabular math columns
\usepackage{nomencl}        % nomenclature generation via makeindex
\usepackage{fancyvrb}       % extended verbatim environments
\usepackage{lettrine}       % dropped capital letter at beginning of paragraph
\usepackage[dvips]{dropping}% alternative dropped capital package
\newcolumntype{d}{D{.}{.}{-1}}
\fvset{fontsize=\footnotesize,xleftmargin=2em}

\makeglossary

 \title{A Residual Smoothing Strategy for Accelerating Newton Method Continuation}

 \author{Dimitri J. Mavriplis
         \thanks{Professor; email: mavriplis@uwyo.edu} \\
         {\normalsize\itshape 
         Department of Mechanical Engineering,
         University of Wyoming, Laramie, WY 82071}
 }

\begin{document}

\maketitle

%------------------------------------------------------------------------------
\begin{abstract} \label{abstract}
A technique for accelerating global convergence of pseudo-transient continuation Newton methods is proposed based on residual smoothing.
The technique is motivated by the effectiveness of local nonlinear smoothers at overcoming strong nonlinear transients.
In the limit of a small pseudo-time step, the method reduces to a local nonlinear smoothing technique,
while in the limit of large pseudo-time steps, an exact Newton method is recovered along with its quadratic convergence properties.
The formulation relies on the addition of a smoothing term source term while leaving the Newton Jacobian matrix unchanged,
thus simplifying implementations for existing Newton solvers.  The proposed technique is demonstrated
on a steady-state and an implicit time-dependent computational fluid dynamics problem, showing significant
gains in overall solution efficiency.

\end{abstract}

%------------------------------------------------------------------------------
\section {Introduction}

Newton methods have become popular strategies for solving large-systems of non-linear equations.
In the field of computational fluid dynamics (CFD) for aerodynamics, Newton methods have enjoyed a resurgence in popularity,
largely due to their ability to provide deep convergence levels for stiff systems of equations
such as those resulting from emerging continuous and discontinuous Galerkin discretizations
particularly for highly resolved steady-state problems  \cite{diosady2009preconditioning,anderson:2017,glasby:2017,ahrabi:2018}.
In the final stages of convergence, when the iterative solution state is close to the exact nonlinear solution,
Newton methods converge in a small number of nonlinear steps, and each step can often be solved effectively
using preconditioned Krylov methods which are generally robust \cite{knoll:2004,mavriplis:nsu3d:2014,boeing:2013}.
However, for most cases, a continuation strategy must be employed to iteratively approach the nonlinear state
where fast convergence is obtained.
A typical strategy consists of employing a pseudo-transient approach
where a pseudo-time term is added to the diagonal of the Jacobian matrix with a variable pseudo-time step \cite{kelley:1998,ceze:2015,ahrabi:2018}.
In the initial phases of convergence, when the pseudo-time step is small, the method approximates
a pseudo-time explicit scheme, and in the final stages, when the pseudo-step becomes large, the exact Newton method is recovered.

Although pseudo-transient continuation for Newton methods is widely used for CFD type problems,
experience has shown that this approach can be slow to converge and often evolves through nonlinear states
that produce convergence stagnation or even unrealizable solutions.
A notable mode of failure or inefficiency occurs when isolated residuals in the nonlinear system remain large
and retard the entire global solution.  This has prompted research into methods such as nonlinear
elimination \cite{lanzkron:1996,cai:2011}, or nonlinear preconditioning \cite{cai:2002,hwang:2005,skogestad:2013,liu:2015}
which attempt to break up the problem into smaller more local nonlinear problems,
based on the observation that local nonlinear solution methods can often overcome the difficulties encountered
by global continuation Newton methods.

Alternatively, the ability of localized residuals to retard the entire global nonlinear problem
may be attributed to a non-smooth residual distribution, since for a perfectly smooth residual distribution,
it should be possible to advance the global problem uniformly. In this work, we propose a residual
smoothing technique that is applied to the right-hand side (i.e. residual vector) of a traditional
pseudo-transient continuation Newton-Krylov method.  Residual smoothing can be achieved through
a variety of approaches, such as using well-known local nonlinear smoothers including block-Jacobi and line smoothers.
In the limit of small pseudo-time steps, the nonlinear iterations of the smoothed Newton continuation
scheme correspond to a nonlinear implicit solver, while in the limit of large pseudo-time steps,
the exact Newton scheme is recovered, along with its quadratic convergence properties.
In the following section, we first begin with the motivation for the current formulation which is given in section \ref{sec:2.1}.
In Section \ref{sec:2.2}, this formulation is generalized and two possible
interpretations are discussed. In Section \ref{sec:3} example
solutions are given for a steady-state and a time-dependent aerodynamic  Reynolds-averaged Navier-Stokes (RANS) CFD
problem on unstructured meshes, demonstrating a more robust and efficient solution strategy compared
to the traditional pseudo-transient continuation Newton-Krylov approach \cite{kelley:1998,ceze:2015,ahrabi:2018}.
In the conclusion section, prospects for further improvements to this approach are discussed.

\section{Formulation}
\subsection{Motivation of the approach}\label{sec:2.1}

We are interested in solving a non-linear set of equations denoted as
\begin{eqnarray}
R(w) = 0
\label{eqn:1}
\end{eqnarray}
where R represents the residual vector which is a nonlinear function of the state vector w for which we seek the solution.
In this work, the residual $R(w)$ may arise from the spatial discretization of a steady-state CFD problem,
or from the equations to be solved at each step of an implicit time-dependent CFD problem.
The Newton scheme for this problem corresponds to solving multiple linear problems of the form
\begin{eqnarray}
\left[ {\partial R(w^n) \over \partial w^n}  \right] \Delta w^n = -R(w^n) 
\label{eqn:2}
\end{eqnarray}
with nonlinear updates given as
\begin{eqnarray}
w^{n+1} = w^n + \alpha \Delta w^n     \,\,\,\,\,\,\,\, with \,\,\,    0 \leq \alpha  \leq 1
\label{eqn:3}
\end{eqnarray}
where the parameter $\alpha$ is determined by a line search which seeks to minimize
the L2 norm of the residual vector $R(w^{n+1})$.
In most cases, a continuation method is required in order to advance the solution through initial nonlinear transients
in order to reach a state close enough to the final solution where the quadratic convergence of Newton's method is observed. 
Pseudo-transient continuation (PTC) is an often used approach which seeks to mimic the physical time evolution
of the solution by adding a pseudo-time term of the form:
\begin{eqnarray}
{M \over \Delta \tau} ( w^{n+1} - w^n) + R(w^{n+1}) = R_t(w^{n+1}) = 0
\label{eqn:4}
\end{eqnarray}
Here $M$ denotes a suitable mass matrix, $\Delta \tau$ represents the pseudo-time step
and $R_t$ corresponds to the augmented pseudo time-dependent residual.
Since the objective is to obtain the solution of $R(w) = 0$, pseudo-time accuracy is not a concern
and a simple first-order pseudo-time discretization (BDF1) is suitable. 
Newton's method is applied to equation (\ref{eqn:4}) in the standard manner, by linearizing $R_t(w)$ about the current
state and solving for the update as:
\begin{eqnarray}
\left[ {M \over \Delta \tau} + \left[ {\partial R(w^n) \over \partial w^n} \right] \right]  \Delta w^n= - R_t(w^{n})
\label{eqn:5}
\end{eqnarray}
Here we note that the pseudo-unsteady residual evaluated at the current state $R_t(w^{n})$ is equal to the original residual
and the above equation can be rewritten as:
\begin{eqnarray}
\left[ {M \over \Delta \tau} + \left[ {\partial R(w^n) \over \partial w^n} \right] \right]  \Delta w^n= - R(w^{n})
\label{eqn:6}
\end{eqnarray}
Finally, the nonlinear update proceeds as determined previously by equation (\ref{eqn:3}) where $\alpha$ can
be determined by a line search.  In this case, the line search seeks to minimize the L2 norm of the pseudo-unsteady
residual, i.e. $F(\alpha) = || R_t(w + \alpha \Delta w) ||_2$, with
\begin{eqnarray}
R_t(w + \alpha \Delta w) = {M \over \Delta \tau } \alpha \Delta w^n + R(w^n + \alpha \Delta w^n)
\label{eqn:7}
\end{eqnarray}
Here we emphasize that it is the (pseudo) time-dependent residual $R_t$ that should
be minimized by the line search, instead of the residual $R$ itself. 
In fact, it is well known that the line search is guaranteed to produce a descent direction
for the L2 norm of this residual provided an exact linearization is used and equation (\ref{eqn:6}) is solved exactly.
This can be shown by denoting the Jacobian matrix as $A = \left[ {M \over \Delta \tau} + \left[ {\partial R(w^n) \over \partial w^n} \right] \right]$, 
and premultiplying equation (\ref{eqn:5}) by $\Delta w^T A^T$, obtaining (dropping the n-superscripts for clarity) \cite{ceze:2015}:
\begin{eqnarray}
\Delta w^T A^T A \Delta w = - \Delta w^T A^T R_t
\label{eqn:8}
\end{eqnarray}
Since the left-hand side of the above equation is a positive scalar quantity, it follows that
\begin{eqnarray}
\Delta w^T A^T R_t(w)  = {1 \over 2} \Delta w^T {d ||R_t(w)||_2 \over dw} \approx {1 \over 2} \Delta ||R_t(w)||_2 < 0
%d ||R_t(w)||_2 \over dw = 2 A^t R_t(w) 
\label{eqn:9}
\end{eqnarray}
Therefore, the solution $\Delta w$ of the linear problem arising from Newton's method in equation (\ref{eqn:5}) or (\ref{eqn:6})
is a descent direction for $||R_t(w)||_2$ and a decrease in this residual is guaranteed for small enough $\alpha$,
thus ensuring the success of the line search.

In the PTC approach, the pseudo-time step $\Delta \tau$ can either represent a global or local time step.
Local time-stepping is most appealing, since we are not interested in time accuracy and since this approach
should result in faster convergence rates, particularly for problems with large variations in mesh cell sizes as is
typical in the applications considered in this work.
In these cases, a pseudo-time step at each mesh cell i is computed as:
\begin{eqnarray}
\Delta \tau_i = CFL . {\Delta \tau_{{explicit}_i}}
\label{eqn:10}
\end{eqnarray}
where ${\Delta \tau_{{explicit}_i}}$ represents an estimate of the explicit time step limit in each cell
and CFL is a global scalar that controls the magnitude of all local time steps. (i.e. a CFL number).
During the continuation process, the CFL number starts out at a low value, of the order of 1, and grows
at each nonlinear step, such that it reaches a large value at the end of the convergence process,
thus recovering an exact Newton scheme for the steady residual problem.

In practice, growth of the CFL number is controlled based on the results of the line search process.
If a full non-linear update ($\alpha = 1.0$) is taken or if $\alpha > 0.75$,
the CFL value is amplified by the factor $\beta_{CFL1} > 1$,
whereas for small values $\alpha \leq  0.1$,
the update is rejected and the CFL is reduced by the factor $\beta_{CFL2} < 1$.
For intermediate values of $\alpha$ the CFL value remains unchanged.
Based on experience, the values $\beta_{CFL1} = 1.5$ and $\beta_{CFL2} = 0.1$
are generally used, since these have been found to provide a good compromise between
speed of convergence and robustness, although other settings are possible.

By its very construction, this CFL-controlled PTC Newton method can never result in solver divergence.
However, stagnation can occur when the solver gets into a nonlinear state from which it cannot recover,
and the CFL is continuously reduced to very small values ($CFL << 1$) producing vanishingly small nonlinear updates.
As an example, this may occur when the solution approaches a non-physical state which would result in localized
negative pressure or density values in computational fluid dynamics problems.

In order to better understand the behavior of the PTC Newton method, we consider the limiting cases of large
and small $\Delta \tau$ values.  Clearly, for large $\Delta \tau$, equation (\ref{eqn:6}) approaches equation (\ref{eqn:2})
and an exact Newton scheme is recovered. On the other hand, for small $\Delta \tau$, equation (\ref{eqn:6}) becomes
\begin{eqnarray}
\left[ {M \over \Delta \tau}  \right]  \Delta w^n= - R(w^{n})
\label{eqn:11}
\end{eqnarray}
which results in updates of the form
\begin{eqnarray}
\Delta w^n= -  \Delta \tau M^{-1} R(w^{n})
\label{eqn:12}
\end{eqnarray}
For a finite-volume scheme where the mass matrix corresponds to the cell volume, this reduces to an explicit
time-stepping scheme (using local time steps in this case).
Even though this should be (linearly) stable for values of $CFL \simeq 1$,
there are many situations where transient nonlinear states coupled with the PTC controller result in values $CFL << 1$
thus producing solver stagnation.
Surprisingly, for many cases where solver stagnation occurs, experience has shown that
commonly used local nonlinear solvers (i.e. non-linear Jacobi, Gauss-Seidel or line solvers)
converge reliably and without difficulty.  In some sense, equation (\ref{eqn:12}) represents
one of the weakest and most inefficient approaches for converging a nonlinear problem.
For example, when solving coupled systems of PDE's, a point implicit approach is much more effective
than a scalar explicit approach.  A simple point implicit approach can be written as
\begin{eqnarray}
\Delta w^n= -  D^{-1} R(w^{n})
\label{eqn:13}
\end{eqnarray}
where $D^{-1}$ represents the inverse of the diagonal block-Jacobian at each mesh point or cell
coupling all equations together.
The justification for the PTC approach is that the non-linear solver will spend very little time in this regime
and that the CFL value will grow quickly, transitioning the solver to a more implicit approach.
However, experience has shown that stagnation at small CFL values is a common failure mode, and even
for problems that converge reliably, a large portion of the overal computational time is spent in the initial
transient regime at low to moderate CFL values.
Furthermore, we speculate that updates of the form given by equation (\ref{eqn:12}) result in poor error smoothing
properties leading to nonlinear states that produce linear systems for Newton's method that are ill-conditioned
(difficult to solve) as well as isolated residuals that impede global convergence and sometimes lead to stagnation.
On the other hand, local nonlinear solvers can be designed specifically for good error smoothing properties,
as is usually done for nonlinear multigrid problems.  A simple example is provided by the Runge-Kutta time-stepping
schemes devised for FAS multigrid methods in CFD \cite{jam2,jam3}.
Particularly, for problems with highly anisotropic meshes (or matrix coefficients),
explicit methods as well as point-implicit methods are known to have poor smoothing properties
and directional or line solvers are required \cite{brandt84,oosterlee_mg,mav16}.

Although Newton methods provide a reliable approach for attaining deep convergence levels for stiff nonlinear problems,
the issues encountered with these methods for evolving the solution through initial nonlinear transients, coupled
with the success of simpler local nonlinear solvers in this regime, have prompted research into nonlinear
preconditioning methods and other approaches that attempt to break up the problem
into smaller more local nonlinear problems \cite{cai:2002,hwang:2005,skogestad:2013,liu:2015}.
A simpler approach consists of initiating the solution process with a local nonlinear solver, and then switching to
a Newton scheme at some point along the solution process.  Aside from lacking elegance, this approach raises various
questions such as how to determine the optimum point for switching solvers, and whether or when to return to the nonlinear
solver should the Newton solver fail to converge \cite{mavriplis:nsu3d:2014}.

In the context of PTC, one would like to design a solver that behaves as a local nonlinear solver rather than
an explicit scheme in the limit of small $\Delta \tau$, while still recovering an exact Newton solver in the limit of large $\Delta \tau$.
Using equation (\ref{eqn:13}) to denote a generic local nonlinear solver, where $D$ may represent a reduced
Jacobian (i.e. block diagonal, lower triangular or tridiagonal line structure)
or more generally any preconditioner that approximates the Jacobian,
one approach may be to modify or replace the mass matrix in equation (\ref{eqn:6}) as:
\begin{eqnarray}
\left[ { D  \over \Delta \tau} + \left[ {\partial R(w^n) \over \partial w^n} \right] \right]  \Delta w^n= - R(w^{n})
\label{eqn:14}
\end{eqnarray}
However, in the limit of small $\Delta \tau$, this results in updates of the form
\begin{eqnarray}
\Delta w^n= -  D^{-1} \Delta \tau R(w^{n})
\label{eqn:15}
\end{eqnarray}
Unfortunately, the updates are still proportional to the pseudo-time step $\Delta \tau$, with the result that stagnation may still occur.
In essence, what is needed is a formulation that reduces to equation (\ref{eqn:13}) in the limit of small $\Delta \tau$
while still recovering Newton's method for large $\Delta \tau$.
This can be achieved using a scheme of the form:
\begin{eqnarray}
\left[ \alpha(\Delta \tau) { D } +  \beta(\Delta \tau) \left[ {\partial R(w^n) \over \partial w^n} \right] \right]  \Delta w^n= - R(w^{n})
\label{eqn:16}
\end{eqnarray}
where $ \alpha(\Delta \tau)$ and $\beta(\Delta \tau)$ are constructed as :
\begin{eqnarray}
 \alpha(\Delta \tau) = {1 \over 1 + \Delta \tau}                                  \\
 \beta (\Delta \tau) = {\Delta \tau \over 1 + \Delta \tau }
\label{eqn:18}
\end{eqnarray}
in order to have the desired asymptotic behavior
\begin{eqnarray}
\alpha(\Delta \tau) \rightarrow 1 \,\,\,\,\, \beta(\Delta \tau) \rightarrow 0 \,\,\,\, for  \,\, \Delta \tau << 1 \\
\alpha(\Delta \tau) \rightarrow 0 \,\,\,\,\, \beta(\Delta \tau) \rightarrow 1 \,\,\,\, for  \,\, \Delta \tau >> 1 
\label{eqn:20}
\end{eqnarray}

There are two undesirable issues with the above formulation.
Firstly, recalling that $D$ may represent any suitable preconditioning matrix,
this approach results in a modified left-hand-side matrix,
the properties of which may be different than the original Jacobian matrix,
particularly in regions $\Delta \tau \sim 1$.
This may require a redesign of the linear iterative solution techniques used to
(approximately) invert the Jacobian, which constitutes a drawback from an implementation viewpoint.
Secondly, and perhaps more importantly, the left-hand side Jacobian in equation (\ref{eqn:16}) no longer
corresponds to an exact linearization of the right-hand side, with the result that the line search procedure
cannot be guaranteed to produce a descent direction for the residual $R(w^n)$.

An alternative approach consists of retaining the same left-hand side as in equation (\ref{eqn:6}) but modifying the right-hand side as:
\begin{eqnarray}
\left[ {M \over \Delta \tau} + \left[ {\partial R(w^n) \over \partial w^n} \right] \right]  \Delta w^n= - R(w^{n}) - D^{-1} {M \over \Delta \tau}  R(w^{n})
\label{eqn:21}
\end{eqnarray}
Clearly, with this formulation, when $\Delta \tau << 1$, the local nonlinear update form of equation (\ref{eqn:13}) is recovered,
while in the limit $\Delta \tau >> 1$, the exact Newton scheme is obtained.
The fact that this formulation can be interpreted as a residual smoothing approach can be seen by re-writing equation (\ref{eqn:21}) as:
\begin{eqnarray}
\left[ {M \over \Delta \tau} + \left[ {\partial R(w^n) \over \partial w^n} \right] \right]  \Delta w^n= - \left[ I + D^{-1} {M \over \Delta \tau}  \right] R(w^{n})
\label{eqn:22}
\end{eqnarray}
and noting that the term $D^{-1} {M \over \Delta \tau}$ is non-dimensional, since both $D$ and ${M \over \Delta \tau}$
scale as the Jacobian matrix. Therefore, the right-hand side in equation (\ref{eqn:22}) is seen to be an average
of residuals with the averaging determined by the sparsity pattern of $D^{-1}$.
This approach has the advantage that the left-hand side Jacobian of the system remains identical to that used in the original
PTC continuation Newton method, with the result that the same linear solver techniques can be used without modification.
Rather the only change is the addition of a constant source term on the right hand side constructed as the scaled corrections
resulting from a local nonlinear solution step. Therefore, given an existing Newton approach,
the implementation of this scheme is relatively straight-forward.

Additionally, the line search procedure is still guaranteed to yield a descent direction for the relevant
smoothed residual vector $R_{sm}$ given as:
\begin{eqnarray}
R_{sm}(w + \alpha \Delta w) = {M \over \Delta \tau } \alpha \Delta w^n + R(w^n + \alpha \Delta w^n) + D^{-1} {M \over \Delta \tau }  R(w^n)
\label{eqn:23}
\end{eqnarray}
Here it is important to note that the smoothing term which corresponds to the scaled nonlinear updates
is evaluated at the start of the Newton step and held constant throughout the line search.
Therefore, this term drops out in the linearization process and the left-hand side Jacobian remains
an exact linearization of the smoothed residual vector defined above, thus guaranteeing that the solution of the
Newton step corresponds to a descent direction for $R_{sm}$.

\subsection{Generalization and Interpretation}\label{sec:2.2}

In the above discussion, the residual smoothing term proportional to $D^{-1} R(w^n)$
was presented as the result of a local nonlinear update, where $D$ represents a suitable
preconditioning matrix which can be either block diagonal (block-Jacobi), lower triangular (Gauss-Seidel),
tridiagonal (line-solver) or any suitable approximation of the full Jacobian $\partial R(w) \over \partial w$
which is known to have good smoothing properties. In practice, the additional source term is formed
by first computing an update produced by a local nonlinear solver as determined by equation (\ref{eqn:13}),
multiplying by $M \over \Delta \tau$ and adding this term to the unsmoothed residual $R(w^n)$.
Under this description, the nature of this smoothing source term can be broadened to include any sequence
of nonlinear updates applied to the original system of nonlinear equations $R(w) = 0$.
For example, the operator $D^{-1}$ may be designed as a sequence of preconditioned nonlinear Runge-Kutta (RK) stages
designed for smoothness following previous work \cite{jam3,mav16,mav17} as:
\begin{eqnarray} \nonumber
w^0 &=& w^n                             \\ \nonumber
w^1 &=& w^0 - \alpha_1 [P]^{-1} R(w^0)    \\ \nonumber
w^2 &=& w^0 - \alpha_2 [P]^{-1} R(w^1)    \\ \nonumber
w^3 &=& w^0 - \alpha_3 [P]^{-1} R(w^2)    \\ 
...                                     \\ \nonumber
w^k &=& w^0 - \alpha_k [P]^{-1} R(w^{k-1})  \\ \nonumber
w^{n+1} &=& w^k                              \nonumber
\label{eqn:24}
\end{eqnarray}
for a k-stage scheme where $\alpha_k$ represent the R-K scheme coefficients
and [P] represents a local preconditioning matrix.
In this case, the final nonlinear update is given as
\begin{eqnarray}
\Delta w^n = w^{n+1} - w^n = w^k - w^0 = -D^{-1}(R(w^n))
\label{eqn:25}
\end{eqnarray}
which implicitly defines the composite nonlinear operator $D^{-1}$.
In this example, the RK scheme corresponds to a sequence of nonlinear updates,
and the matrix [P] may be held fixed or evolved as a function of the nonlinear state $w^k$ at each state.
In practice, any number of nonlinear iterations may be used to form the smoothing source term,
and any nonlinear solver may be used, including multiple passes of preconditioned RK,
nonlinear block-Jacobi, Gauss-Seidel, line solvers,
or even linear or non-linear (FAS) multigrid.
Presumably, more passes will provide smoother residual terms and accelerated non-linear convergence overall.

The proposed solution strategy may be interpreted in two different manners.
On the one hand, it can be seen as a more elegant strategy (compared to switching solvers)
for smoothly transitioning from local nonlinear smoothers, which are well known to provide good initial convergence,
to an increasingly exact Newton method, which is generally more successful in obtaining deep final convergence,
as the pseudo-time step is increased.

The second interpretation is that of residual smoothing,
where the process ensures a smooth field of residuals on the right-hand-side of equation (\ref{eqn:22}),
thus producing smooth updates to the solution vector at each Newton step,
thereby avoiding cases where isolated residuals hold up global convergence
and avoiding evolution of the solution state into unphysical nonsmooth states that either cause failure,
are difficult or costly to recover from using the typical PTC controller strategy of lowering $\Delta \tau$,
or which may drive the pseudo-time step to zero and cause stagnation.
Indeed, using the current line-search based controller,
if the state does become ill-conditioned causing the pseudo-time step to be reduced,
the scheme reverts to the local nonlinear smoother $\Delta w^n = - D^{-1}(R(w^n))$
even in the limit $\Delta \tau << 1$, which results in a lower bound on the degree
to which corrections are reduced with smaller pseudo-time steps, and which is inherently self-correcting.

\section{Results}\label{sec:3}

The proposed residual smoothing approach is illustrated using two computational fluid dynamics (CFD) test cases.
We seek solutions of the Reynolds-averaged Navier-Stokes (RANS) equations, which consist of the compressible
form of the Navier-Stokes equations augmented with the one-equation Spalart-Allmaras turbulence model \cite{spalart,allmaras_iccfd7}.
In three-dimensions, this yields a system of 6 coupled PDE's describing conservation of mass (1), momentum (3),
energy (1), and turbulence convection/diffusion/production (1). These equations are discretized using a vertex-based finite-volume
approach on unstructured meshes of mixed element types, using isotropic tetrahedral elements in off-body regions and highly anisotropic
prismatic elements in near-body regions to capture thin boundary layer gradients, as is typical for external aerodynamics CFD problems.
The finite-volume discretization employs a matrix-based artificial dissipation formulation
which is second-order accurate and which employs a distance-2 (neighbor of neighbor) stencil \cite{mav24,mavriplis:nsu3d:2014}.
This results in a relatively large bandwidth Jacobian matrix which is generally
considered impractical to store in memory. However, exact Jacobian-vector products are available through
hand-differentiation of the residual routine and are used in the Newton-Krylov solver.
Alternatively, the Jacobian of the corresponding first-order accurate discretization has a nearest-neighbor stencil
and may be stored for use as a preconditioner or for direct use in the local nonlinear solvers as described below.

\begin{figure}
 \centering
  \subfigure{\hspace{-0.25in}
  \includegraphics[width=2.0in]{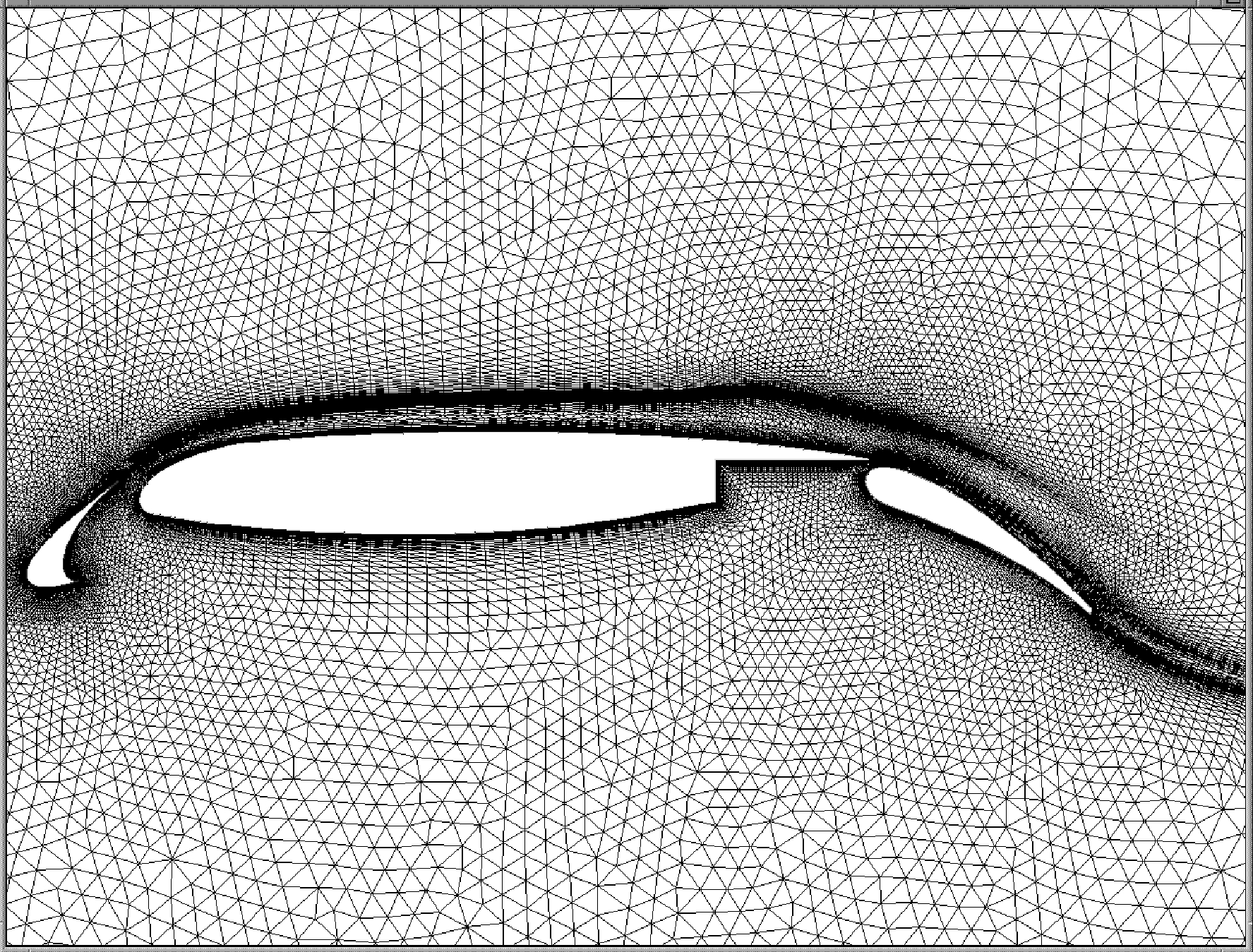}}
  \subfigure{\hspace{-0.00in}
  \includegraphics[width=2.0in]{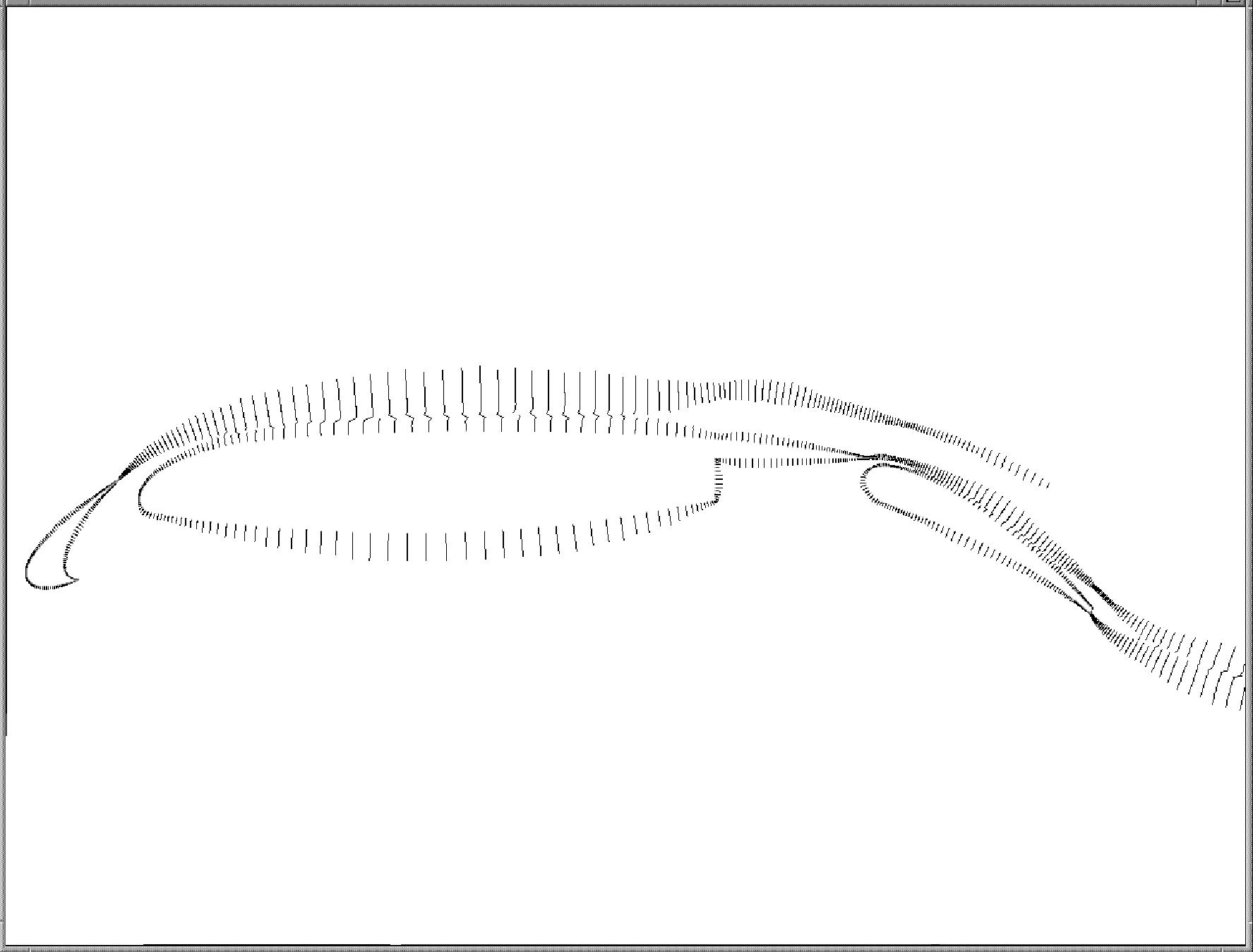}}
  \subfigure{\hspace{-0.25in}
  \raisebox{-1.75in} {\includegraphics[width=2.65in]{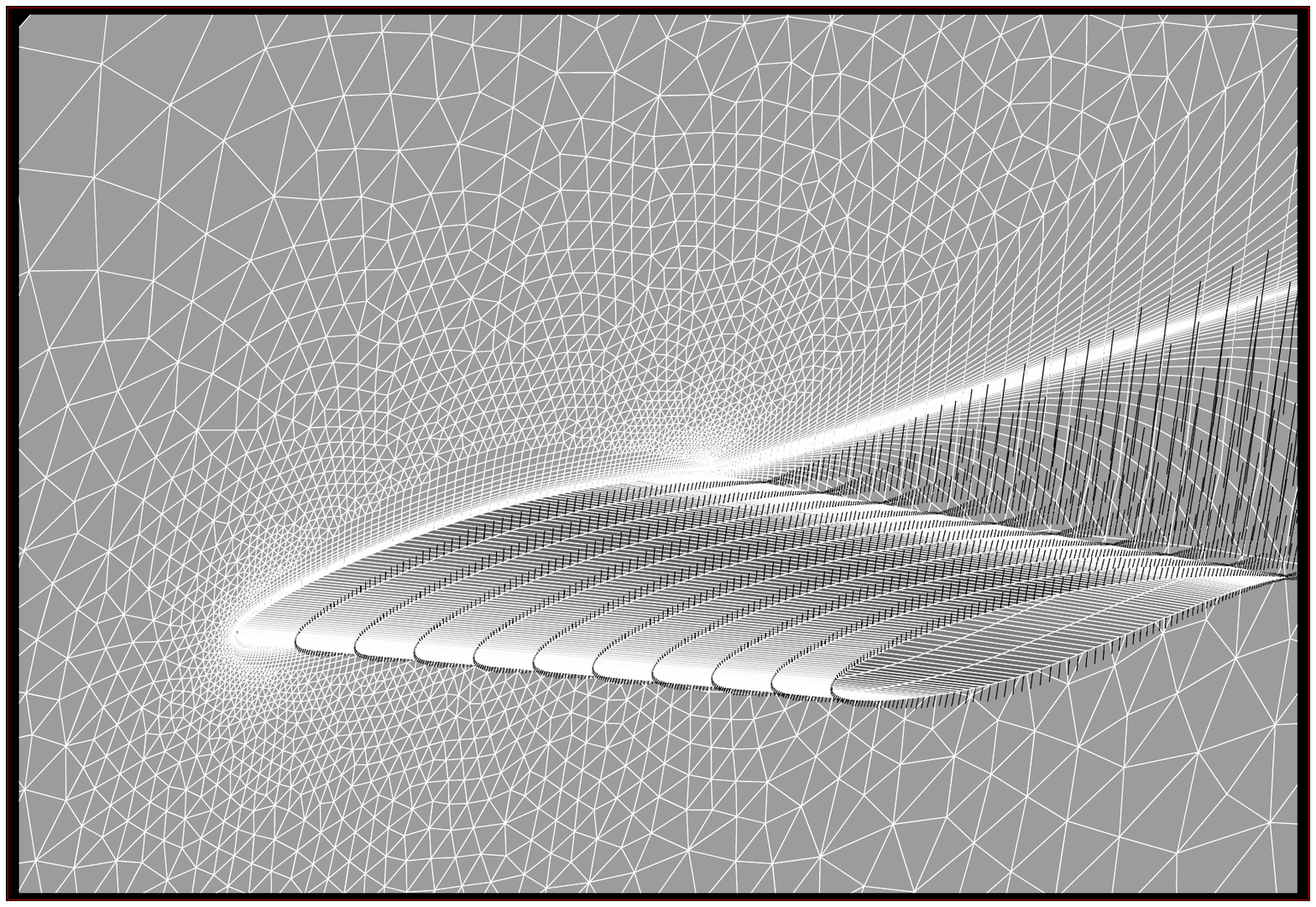}}}
\vspace{-1.75in}
  \caption{Illustration of implicit line structures extracted from unstructured meshes in two and three dimensions}
  \label{fig:lines}
\end{figure}

The baseline local nonlinear solver consists of a 3-stage line-preconditioned Runge-Kutta scheme designed for
good error smoothing properties \cite{mav16,mav17}.
This corresponds to the nonlinear update scheme depicted in equation (\ref{eqn:24})
where[P] is a piecewise block-tridiagonal matrix, as determined by sets of lines constructed in the mesh using a graph algorithm.
The line structures are used to relieve the stiffness associated with high mesh stretching in near-wall regions.
The lines are constructed based
on an algorithm that is initiated in regions of high mesh anisotropy and proceeds in a greedy fashion towards
regions of lower anisotropy, with lines terminating when nearly isotropic mesh elements are encountered \cite{mav16,mav17,mav20}.
This graph algorithm results in a set of lines of varying length that do not
span the entire mesh, as illustrated in Figure \ref{fig:lines}. For each identified edge joining points ij in the line set, the
corresponding first-order accurate Jacobian off-diagonal block matrix entries $[O_{ij}]$ and  $[O_{ji}]$ are retained in the [P] matrix,
while all other diagonal entries are dropped. The mesh vertices and edges are then reordered
resulting in a local tridiagonal Jacobian matrix structure for each line.
The implementation naturally handles
lines of varying length, where lines of length 0 (1 vertex, zero edges) reduce to a block-diagonal preconditioning approach.
In the current implementation, the local tridiagonal matrices are factorized and frozen for all stages of the
Runge-Kutta solver.
This line-preconditioned RK scheme can either be used directly as a nonlinear solver, or as a smoother for a non-linear FAS multigrid
scheme based on coarse agglomerated meshes \cite{mav16,mav17}.

Alternatively, a pseudo-transient continuation (PTC) Newton-Krylov method can be used to solve the same problem.
In this case, a right-preconditioned GMRES algorithm is used to solve the linear system described by equation (\ref{eqn:6})
using exact Jacobian-vector products. Preconditioning is achieved by approximately inverting the equivalent first-order Jacobian matrix
using a small number of linear multigrid iterations driven on each level by a linear block-tridiagonal line solver
which is analogous to the nonlinear variant described previously.
%%%%%%%%%%%%%%%%%%%%%%%%%%%%%%%%%%%%%%%%%%%%%%%%%%%%%%%%%%%%%%%%%%%%%             
%%%%%%%%%%%%%%%%%%%%%%%%%%%%%%%%%%%%%%%%%%%%%%%%%%%%%%%%%%%%%%%%%%%%%     

\subsection{Steady-State Case}
\begin{figure}
  \centering
  \subfigure[Unstructured mesh with near wall detail]{
  \includegraphics[width=0.40\textwidth]{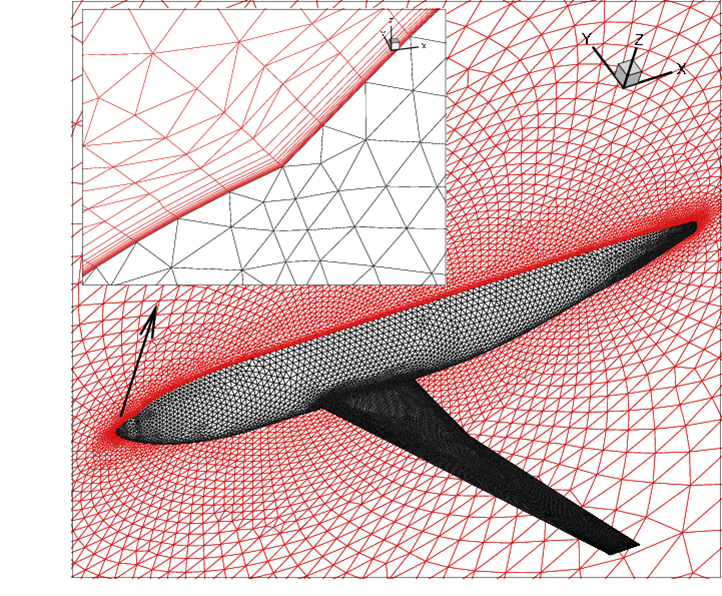} }
  \subfigure[Steady state solution in terms of computed surface pressure contours]{
  \includegraphics[width=0.40\textwidth]{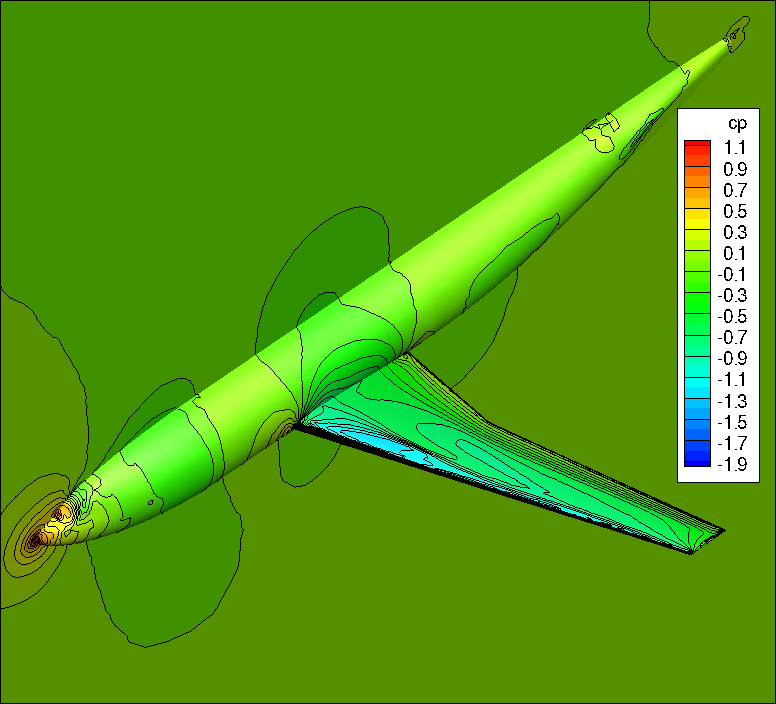} }
  \caption{Illustration of steady-state wing-body CFD test case}
  \label{fig:2}
\end{figure}

The first test case consists of the steady transonic flow over an aircraft wing-body geometry which has been the subject
of previous drag-prediction-workshop (DPW) accuracy studies \cite{dpw3:summary:2008,beth:2009}
as well as convergence efficiency studies \cite{mavriplis:nsu3d:2014}.
An unstructured mesh of 1.2 million vertices is used for this case, as shown in Figure \ref{fig:2}(a) where the highly anisotropic
prismatic elements in the near wall regions are depicted. The aspect ratio of these elements can be of the order of $10^4$.
The freestream flow conditions include a Mach number of 0.75, a flow incidence of 0 degrees,
and a Reynolds number based on the mean aerodynamic chord of the wing of 3 million.
The steady-state solution produced by all the solution techniques is illustrated in Figure \ref{fig:2}(b),
where a weak shock wave on the wing is observed.
This is considered a relatively well behaved (easy to converge) case by industrial standards.

\begin{figure}[!h]
  \centering
  \subfigure[Single grid solver]{
  \includegraphics[width=0.40\textwidth]{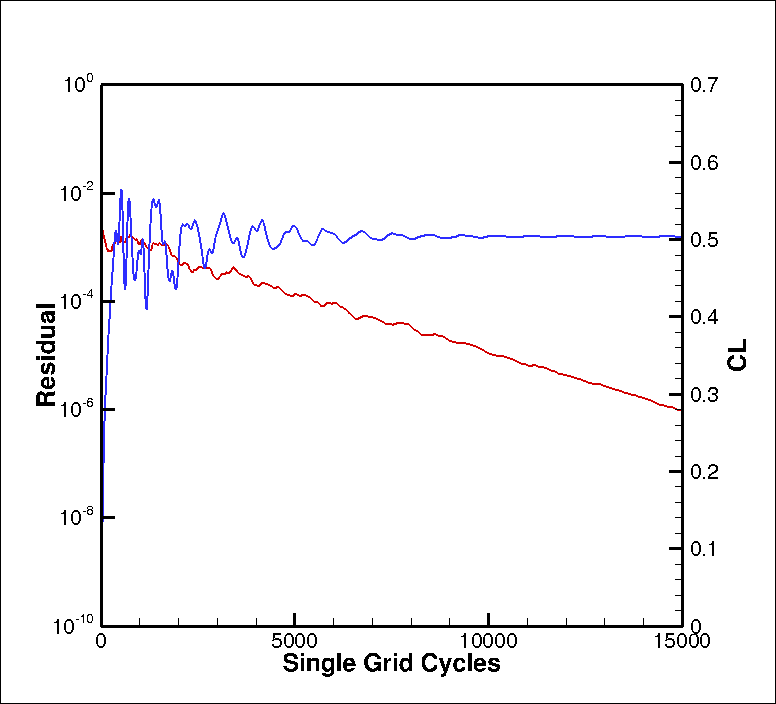} }
  \subfigure[Multigrid solver]{
  \includegraphics[width=0.40\textwidth]{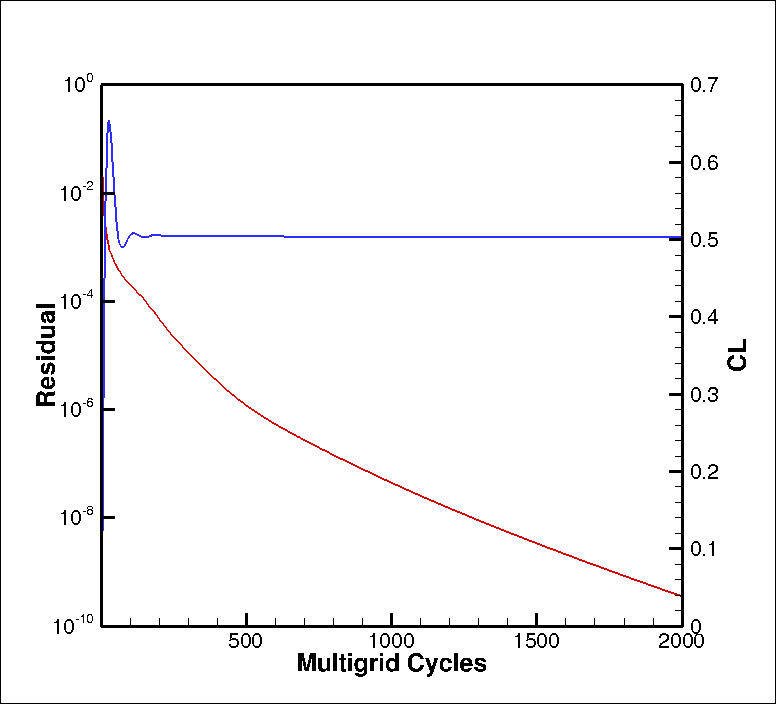} }
  \caption{ Convergence history of (a) single grid and (b) multigrid nonlinear solvers. Note difference in horizontal scales}
  \label{fig:4}
\end{figure}

For all solvers, the flow field is initialized impulsively as a uniform freestream field that does not satisfy the wall boundary conditions.
Figure \ref{fig:4} depicts the convergence history achieved by the nonlinear line-preconditioned RK solver used alone on the fine mesh,
as well as that achieved using this solver as a smoother for the nonlinear multigrid scheme with 4 coarser agglomerated mesh levels.
Convergence is monitored as a function of the L2 norm of the total residual and the history of the computed lift coefficient,
which represents a global integrated quantity of engineering interest.
In both cases, a near monotonic residual convergence is observed, with the multigrid scheme achieving significantly faster
overall convergence, as expected.
%----------------------------------------------------------------------
%----------------------------------------------------------------------

\begin{figure}{!h}
  \centering
  \subfigure[]{
  \includegraphics[width=0.40\textwidth]{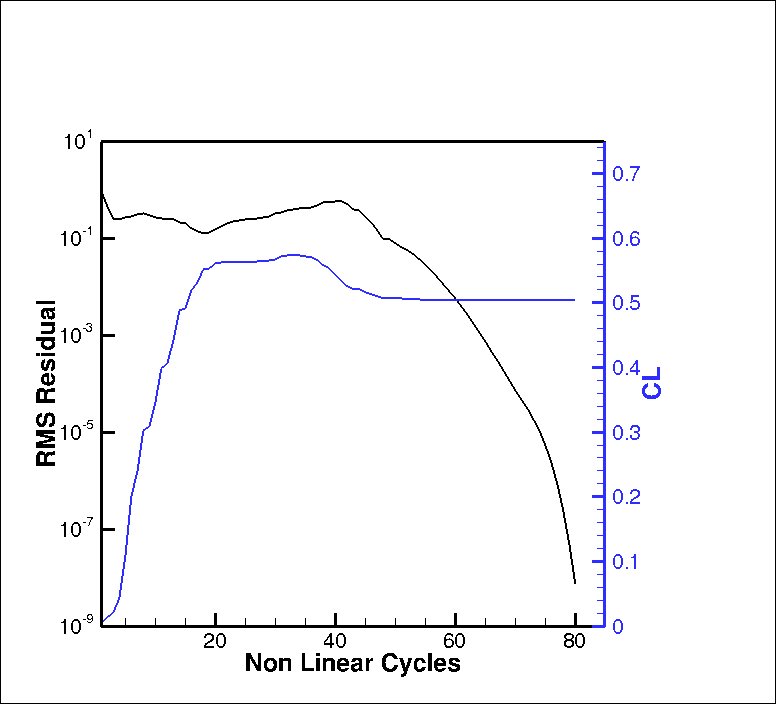} }
  \subfigure[]{
  \includegraphics[width=0.40\textwidth]{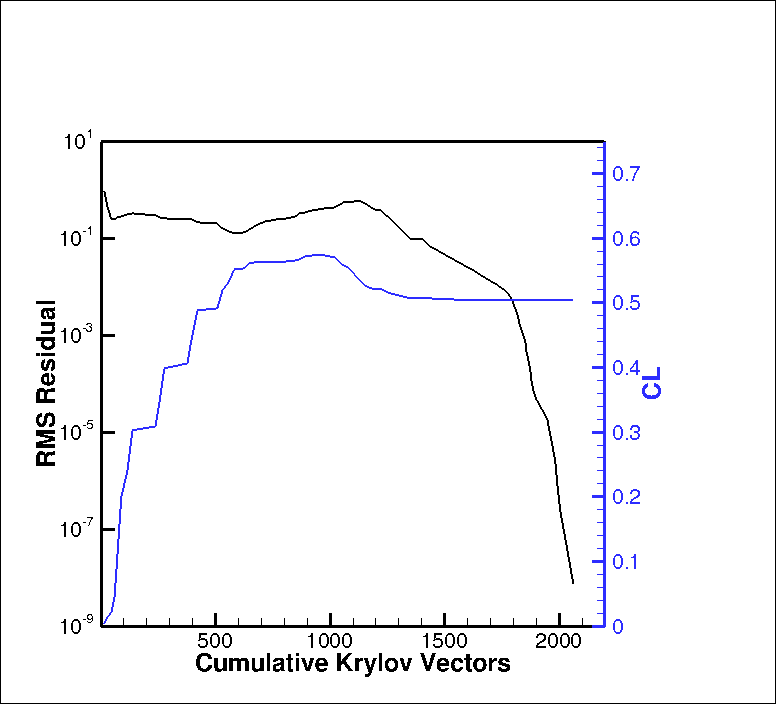} }
  \caption{Convergence history for Newton-Krylov solver (unsmoothed)}
  \label{fig:5}
\end{figure}

Figure \ref{fig:5} depicts the convergence history achieved by the PTC Newton-Krylov solver.
In this case, a linear system convergence tolerance of two orders of magnitude reduction in the linear residual
is requested and a maximum of 100 Krylov vectors are allowed. If the requested convergence tolerance is not met within the maximum
number of Krylov vectors, the linear solver exits and records a failure, which prompts the PTC controller to reject the nonlinear update
and lower the pseudo-time step by the factor $\beta_{CFL2}$. Each Krylov vector is preconditioned using 3 linear multigrid cycles
with 4 line-smoother passes on each level. The computational cost of a Krylov vector is approximately equal to the cost of a single
nonlinear multigrid cycle or two nonlinear single grid cycles in the previous case.
The PTC controller uses an initial pseudo-time step value of CFL=10, and increases or reduces the CFL number depending on the success of the
line search operation based on the amplification factors $\beta_{CFL1} = 1.5$ and  $\beta_{CFL2} = 0.1$, respectively.
From Figure \ref{fig:5}(a), the residual is seen to decrease by 8 orders of magnitude over 80 nonlinear cycles (Newton steps).
However, for the first 40 or more cycles, the residual does not decrease substantially from its initial value,
and only begins to decrease rapidly near the end of the calculation, as is typically observed for Newton methods.
Figure \ref{fig:5}(b) reproduces the same convergence histories plotted in terms of cumulative Krylov vectors,
which is a better representation of computational wall-clock time, further illustrating how the slow initial convergence of the Newton
scheme consumes more than half of the Krylov vectors or compute time.  Additionally, the lift coefficient is also slow to converge initially,
and only comes close to its final value in the quadratic convergence region of the Newton method.
This is in contrast to the rapid initial convergence of this quantity in the nonlinear multigrid solver, as shown in Figure \ref{fig:4},
which may seem surprising given the fact that the Newton-Krylov scheme uses an analogous linear multigrid solver as a preconditioner.
We note that rapid partial convergence of engineering quantities such as integrated force coefficients
can have important practical implications for production use of solvers.
%----------------------------------------------------------------------
%----------------------------------------------------------------------
%----------------------------------------------------------------------
%----------------------------------------------------------------------
%----------------------------------------------------------------------
\begin{figure}[!h]
  \centering
  \subfigure[]{
  \includegraphics[width=0.40\textwidth]{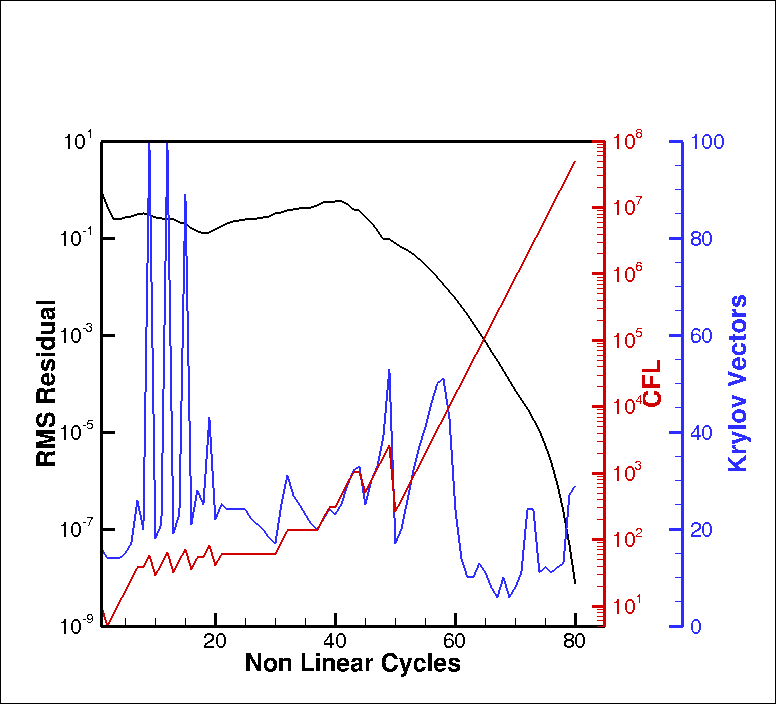} }
  \subfigure[]{
  \includegraphics[width=0.40\textwidth]{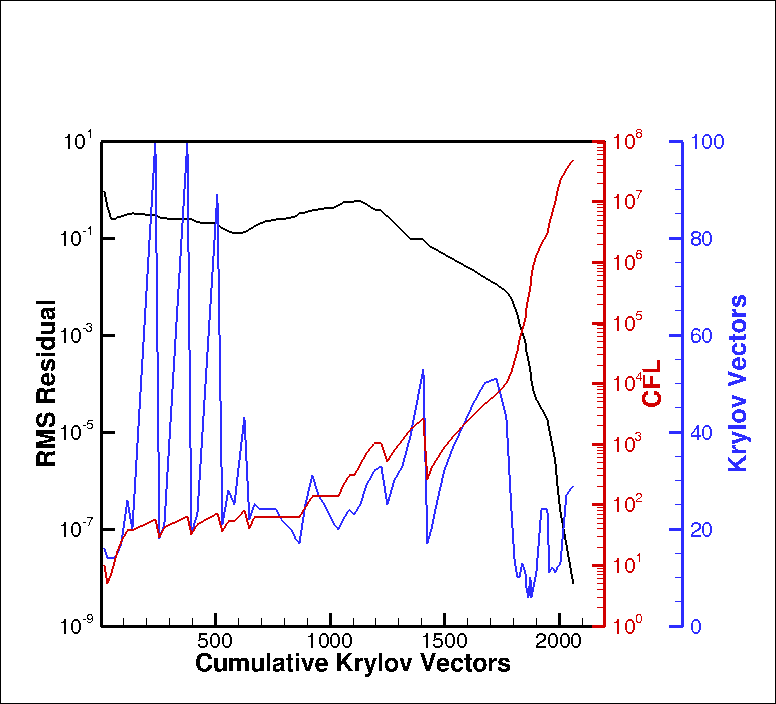} }
  \caption{Convergence history details of pseudo-transient continuation process for Newton-Krylov solver (unsmoothed)}
  \label{fig:6}
\end{figure}
%----------------------------------------------------------------------
%----------------------------------------------------------------------
\begin{figure}[h]
  \centering
  \subfigure[]{
  \includegraphics[width=0.40\textwidth]{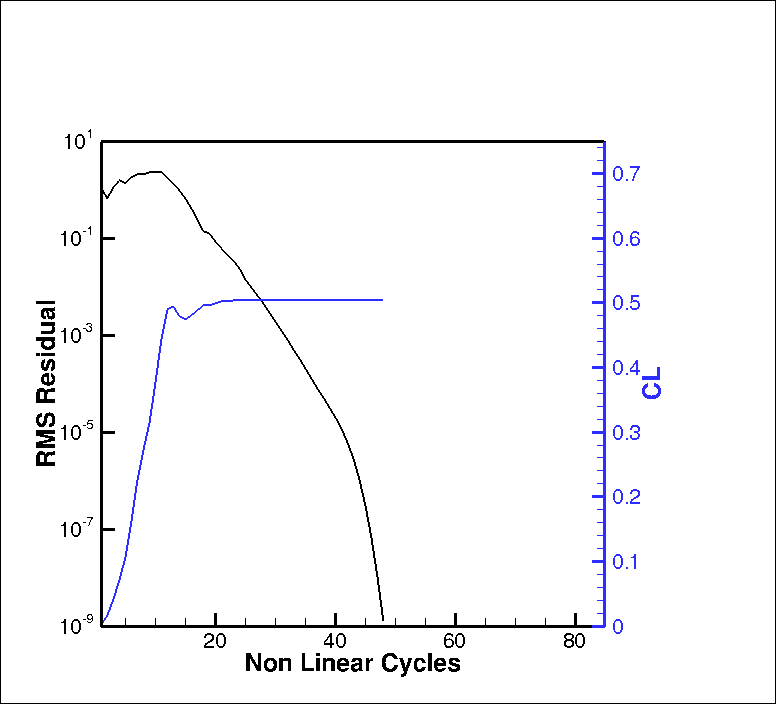} }
  \subfigure[]{
  \includegraphics[width=0.40\textwidth]{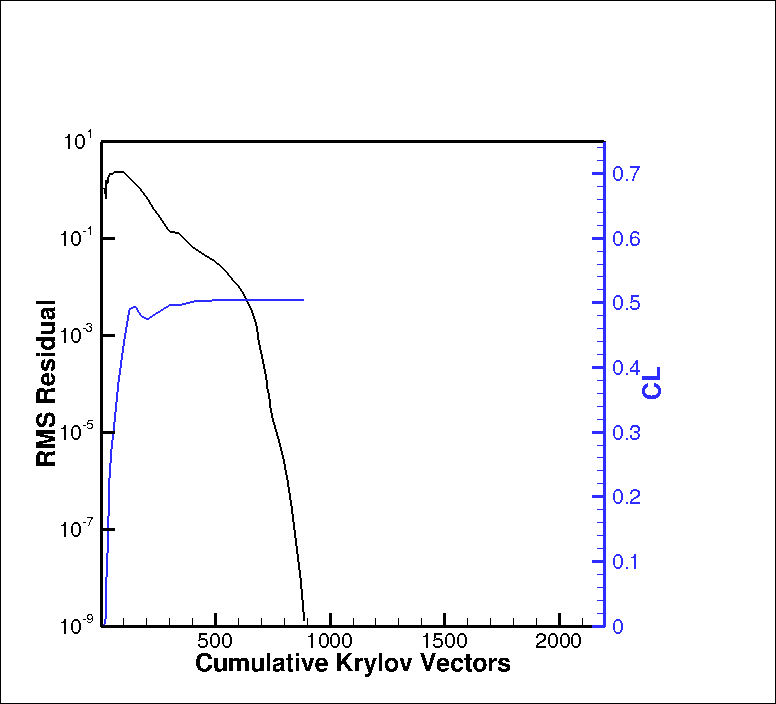} }
  \caption{Convergence history for residual-smoothing Newton-Krylov solver (single grid RK smoother)}
  \label{fig:7}
\end{figure}
%----------------------------------------------------------------------
%----------------------------------------------------------------------
\begin{figure}[!h]
  \centering
  \subfigure[]{
  \includegraphics[width=0.40\textwidth]{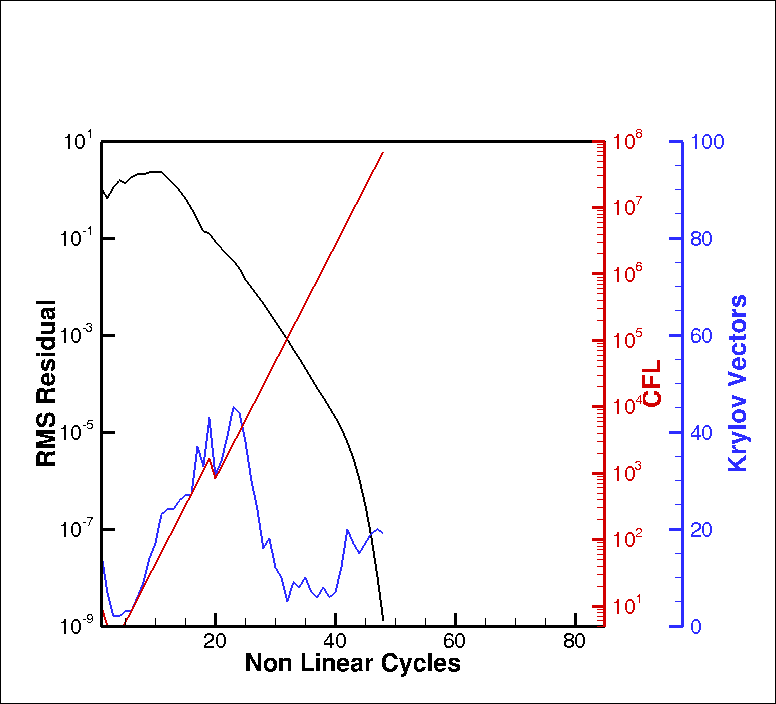} }
  \subfigure[]{
  \includegraphics[width=0.40\textwidth]{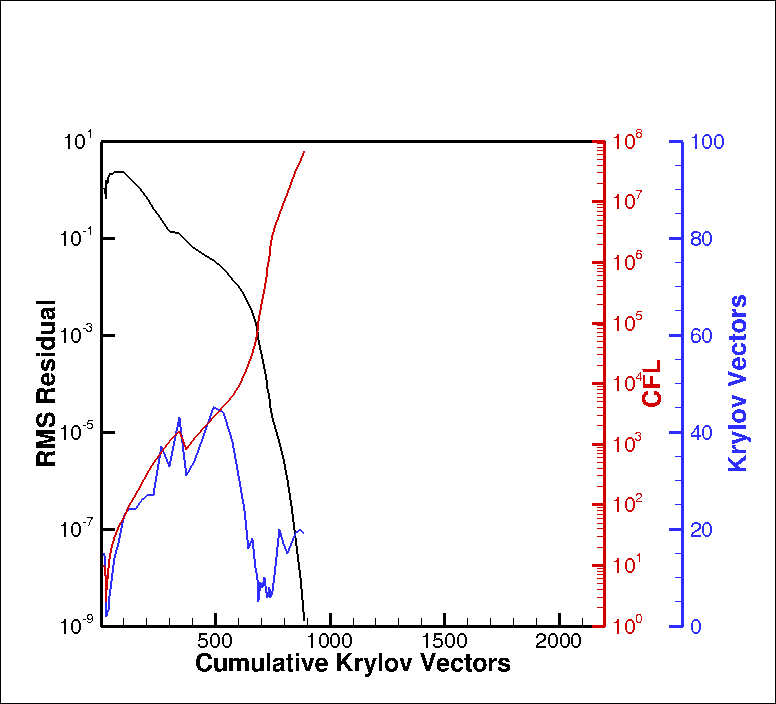} }
  \caption{Convergence history details of pseudo-transient continuation for smoothed Newton-Krylov solver (single grid RK smoother)}
  \label{fig:8}
\end{figure}
%----------------------------------------------------------------------
%----------------------------------------------------------------------
\begin{figure}[!h]
  \centering
  \subfigure[]{
  \includegraphics[width=0.40\textwidth]{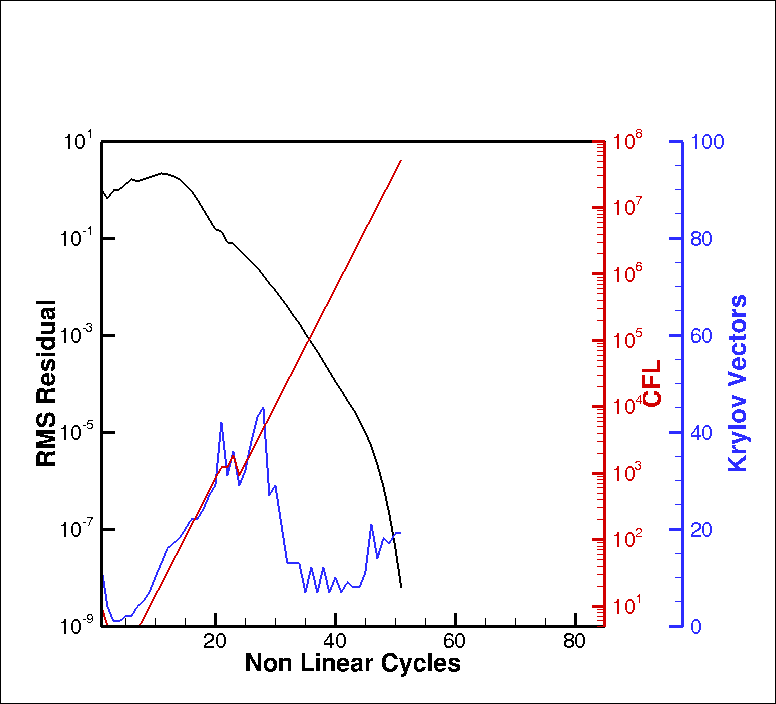} }
  \subfigure[]{
  \includegraphics[width=0.40\textwidth]{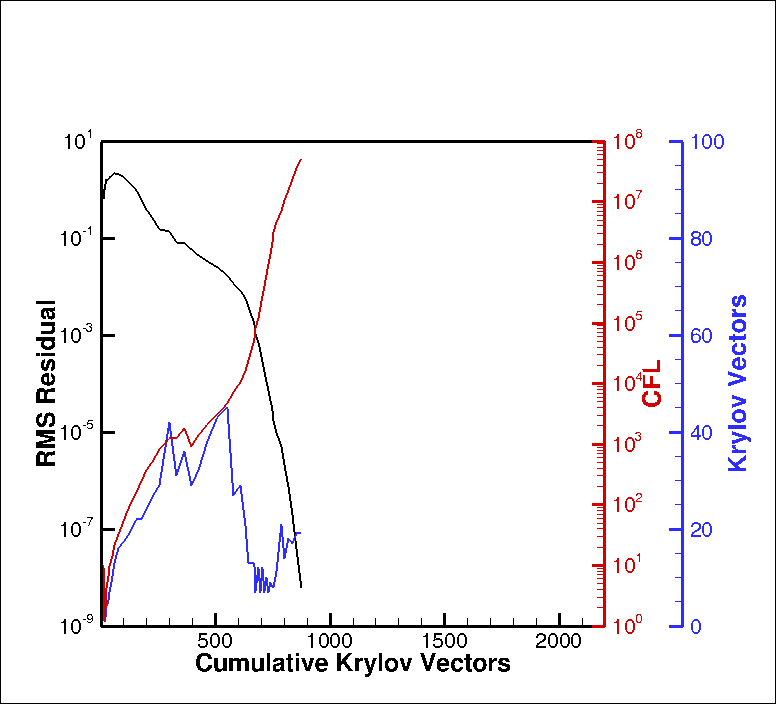} }
  \caption{Convergence history details of pseudo-transient continuation for smoothed Newton-Krylov solver using multigrid smoother}
  \label{fig:9}
\end{figure}
%----------------------------------------------------------------------
%----------------------------------------------------------------------
\begin{figure}[!h]
  \centering
  \subfigure[]{
  \includegraphics[width=0.40\textwidth]{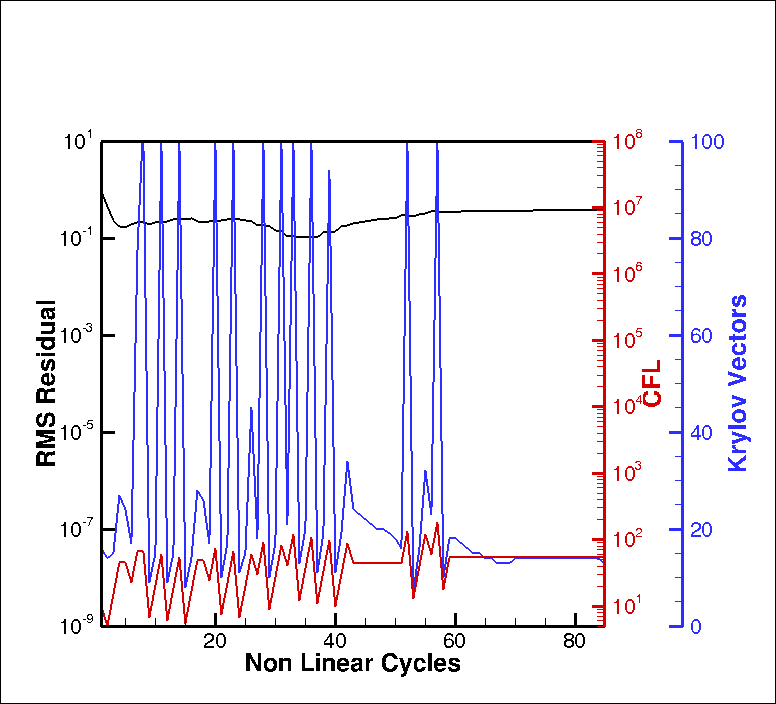} }
  \subfigure[]{
  \includegraphics[width=0.40\textwidth]{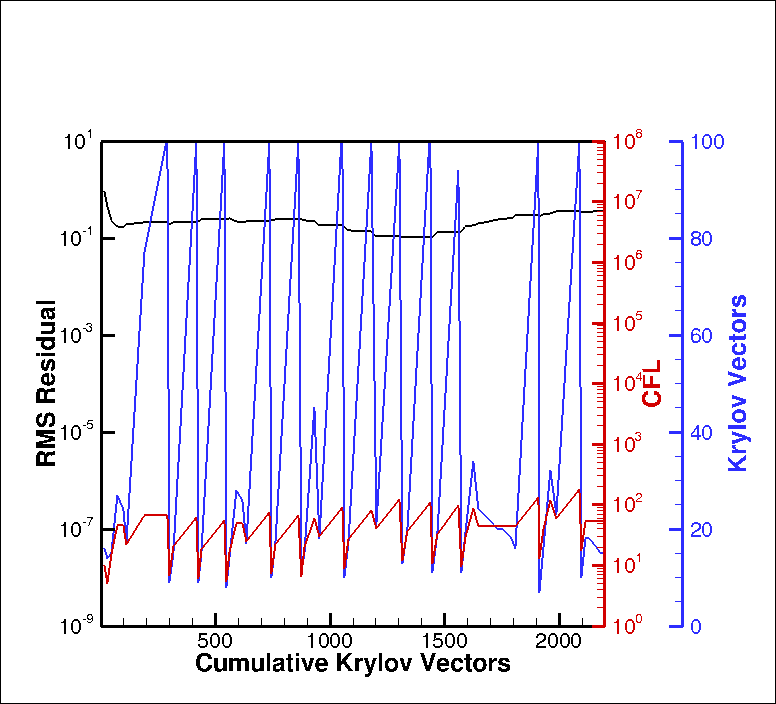} }
  \caption{Convergence history details of pseudo-transient continuation for unsmoothed Newton-Krylov solver using more aggressive CFL growth factor $\beta_{CFL1} = 3$ (single grid RK smoother)}
  \label{fig:10}
\end{figure}
%----------------------------------------------------------------------
%----------------------------------------------------------------------
\begin{figure}[!h]
  \centering
  \subfigure[]{
  \includegraphics[width=0.40\textwidth]{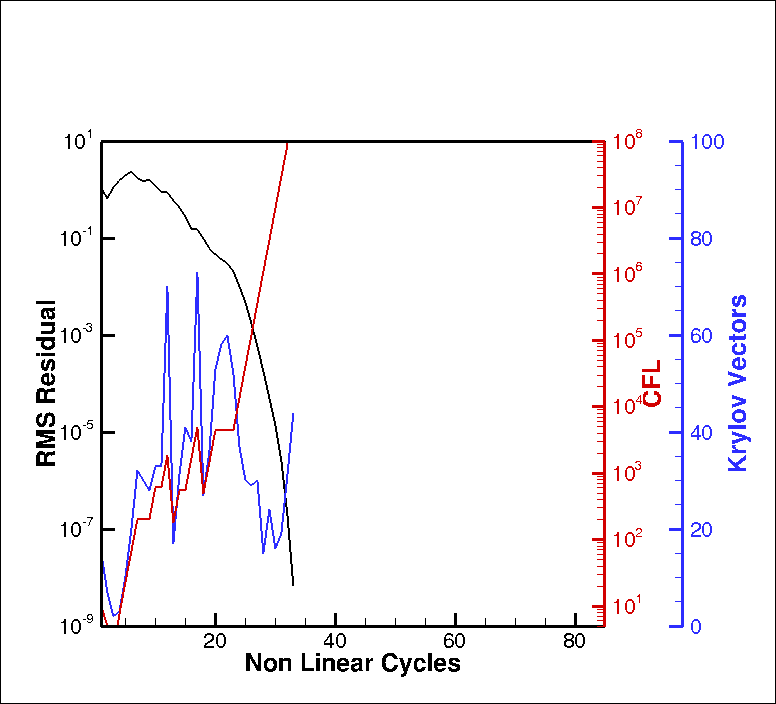} }
  \subfigure[]{
  \includegraphics[width=0.40\textwidth]{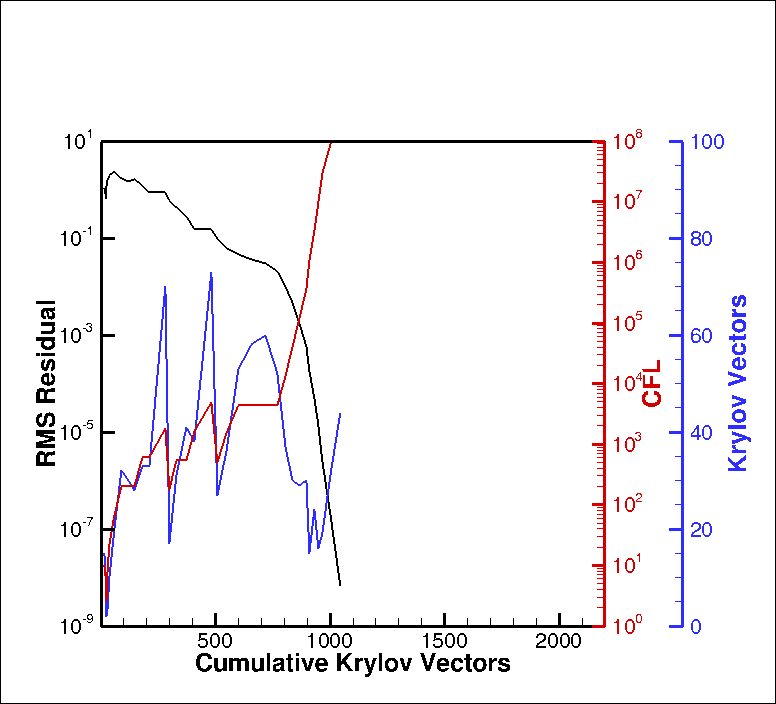} }
  \caption{Convergence history details of pseudo-transient continuation for smoothed Newton-Krylov solver using more aggressive CFL growth factor $\beta_{CFL1} = 3$ (single grid RK smoother)}
  \label{fig:11}
\end{figure}
%----------------------------------------------------------------------

In order to gain a better understanding of the Newton-Krylov convergence process, Figure \ref{fig:6} plots, in addition to the nonlinear residual,
the number of Krylov vectors per Newton step and the CFL value as determined by the PTC controller.
We first note that the initial value CFL=10 is too large for the impulsive initial condition, and is reduced by the PTC controller on the first
step to CFL=1.0. Thereafter the CFL value rises slowly at first,
but stagnates and remains relatively low over the first 50 nonlinear steps (out of 80),
which corresponds to the first 1400 cumulative Krylov vectors (out of 2063), after which it resumes steady growth to large
values that enable the final quadratic convergence. Another point to note is that the hardest linear systems to solve, i.e. those that require
the largest number of Krylov vectors, are those at the beginning of the continuation process, when the CFL or pseudo-time step is small,
and the Jacobian is thought to be strongly diagonally dominant.  Evidently, the nonlinear state about which these Jacobians are linearized
results in stiff linear problems, particularly compared to the final state where quadratic nonlinear convergence occurs. 
In this sense, the pseudo-transient continuation approach is doubly inefficient. On the one hand it is slow to grow the CFL value to produce
a more strongly implicit scheme, while on the other hand it results in linear problems early in the continuation process
that are more difficult and costly to solve than even the final Newton step.

Figure \ref{fig:7} depicts the convergence history achieved by the residual smoothing approach applied to the PTC Newton-Krylov solver.
In this case, 5 nonlinear cycles of the 3-stage line-preconditioned RK scheme are used to generate the smoothing source
term, which is added to the right-hand side as given in equation (\ref{eqn:21}), with all other Newton-Krylov solver parameters remaining
identical to those used in the previous case. As the figure illustrates, the number of nonlinear cycles required
to achieve the same residual reduction of 8 orders of magnitude is reduced from 80 to 48. Perhaps more importantly, the number
of cumulative Krylov vectors is reduced from 2063 to 888. Given that the overhead of the residual smoothing process accounts
for roughly 10\% of the total solver time, the overall gain is roughly a factor of 2.
Furthermore, the initial convergence of the lift coefficient is more rapid,
particularly in terms of the number of cumulative Krylov vectors,
displaying convergence behavior which is closer to that observed in Figure \ref{fig:4} using the nonlinear multigrid approach.
Figure \ref{fig:8} depicts additional details of the solution process for this case.
As in the previous case, the initial value CFL=10 is determined by the line search to be too large for the impulsive initial condition,
and is reduced to CFL=1 at the first step. However, the CFL value grows almost monotonically thereafter, quickly reaching large
values and enabling earlier onset of quadratic nonlinear convergence.  At the same time, the number of Krylov vectors
required to solve each linear system is rarely more than double the required number at the final state, presumably due
to smoother nonlinear states achieved throughout the continuation process.

In Figure \ref{fig:9} the residual smoothing approach is repeated, although the residual smoothing term is
now generated using 5 nonlinear multigrid cycles (using the 3-stage line-preconditioned RK scheme as a smoother on each level).
As demonstrated in Figure \ref{fig:4}, the nonlinear multigrid solver is known to be a much more effective nonlinear solver
than the single grid line-RK smoother alone.
However, from Figure \ref{fig:9} it is seen that the convergence behavior using the multigrid scheme
in the place of the single grid nonlinear solver
for residual smoothing purposes does not improve the overall convergence behavior. In this case, the equivalent
8 order of magnitude residual drop is achieved in 51 nonlinear cycles and 874 cumulative Krylov vectors, compared to 48 and 888 respectively
for the previous case.  This case provides evidence that the effective mechanism of the proposed approach is indeed smoothing
rather than simply advancing the solution nonlinearly through the application of additional nonlinear solver steps,
as would be the case if a solver switching strategy were to be used.

Returning to the original residual smoothing results plotted in Figure \ref{fig:8},
since the CFL grows almost monotonically, it can be argued
that a faster CFL growth rate may lead to faster nonlinear convergence.
In Figures \ref{fig:10} and \ref{fig:11} the original unsmoothed and the smoothed Newton-Krylov solvers
discussed in Figures  \ref{fig:6} and  \ref{fig:8} are rerun using a larger CFL growth factor of $\beta_{CFL1} = 3$.
As seen in Figure \ref{fig:10}, the unsmoothed Newton-Krylov solver fails to converge with this more aggressive
CFL growth factor, entering a limit cycle
leading to repeated reduction in the CFL values, followed by resumed growth.
On the other hand, the smoothed approach achieves convergence (defined as 8 orders of magnitude residual reduction)
in 33 nonlinear cycles corresponding to 1044 cumulative Krylov vectors.
In spite of the fast CFL growth rate and resulting lower number of nonlinear cycles,
the total number of Krylov vectors is slightly higher and the overall efficiency is slightly lower than the convergence discussed
previously and shown in Figures \ref{fig:7} and \ref{fig:8}. On the one hand, this test case demonstrates the ability of the smoothed
approach to overcome convergence stalling that occurs for the nonsmoothed approach, thus providing
an additional level of robustness. However, this case also reveals the delicate balance that
must be achieved between local smoothing and the continuation process in order to maximize overall solution efficiency.

\subsection{Implicit Time Dependent Case}
\begin{figure}[!h]
  \centering
  \subfigure[Far field view of grid]{
  \includegraphics[width=0.30\textwidth]{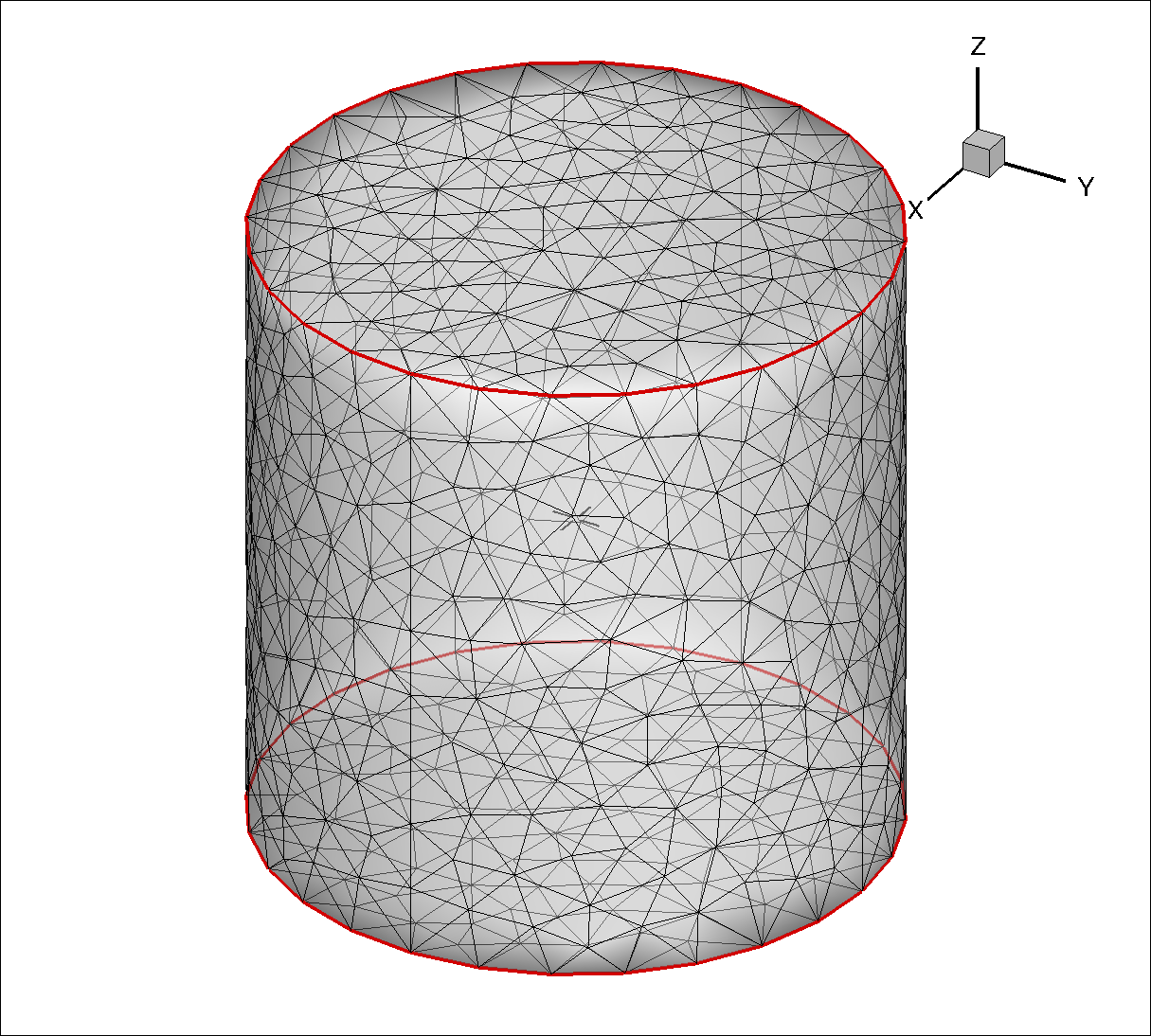} }
  \subfigure[Near field view of grid]{
  \includegraphics[width=0.30\textwidth]{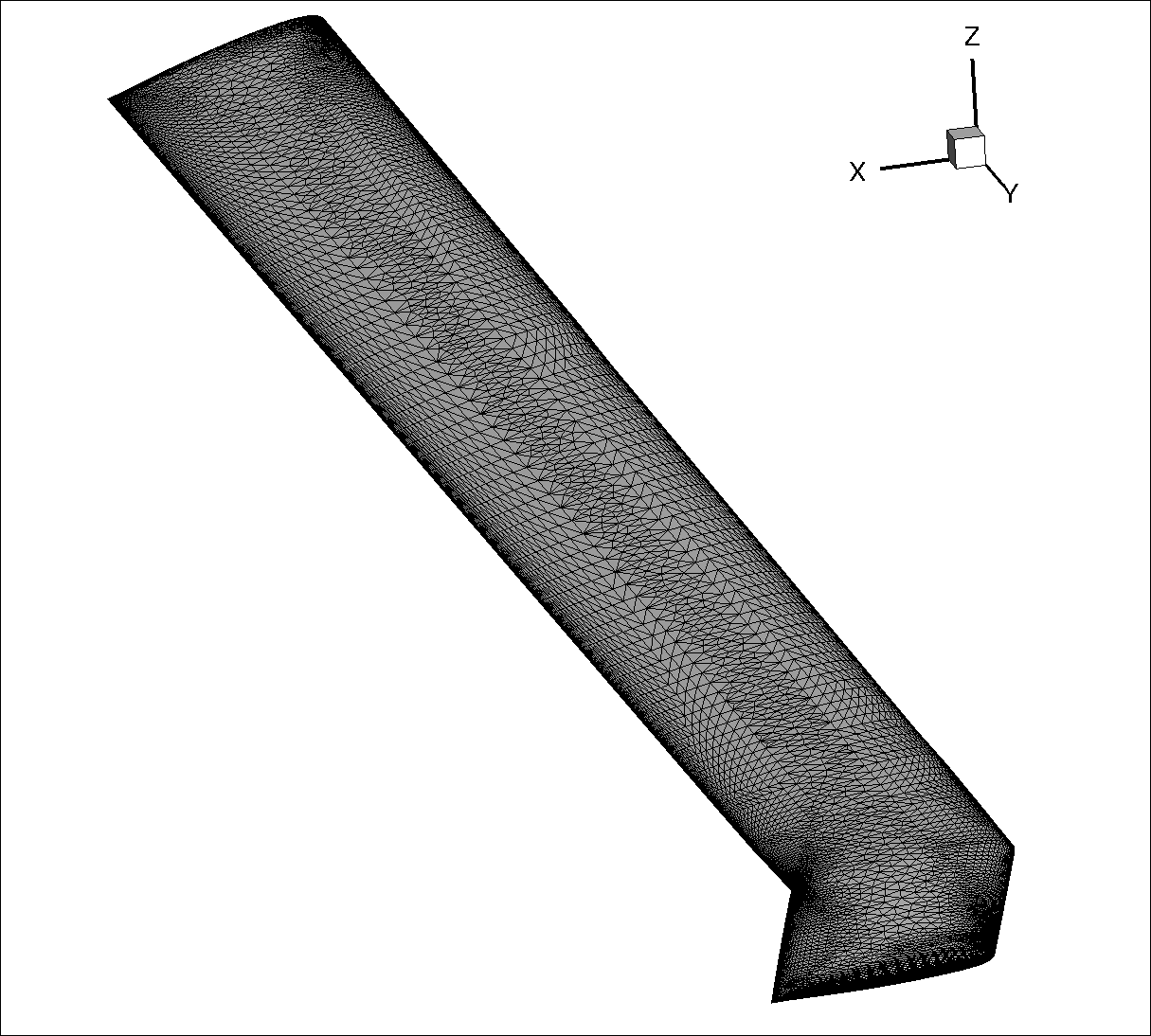} }
  \subfigure[Nonlinear MG convergence]{
  \includegraphics[width=0.30\textwidth]{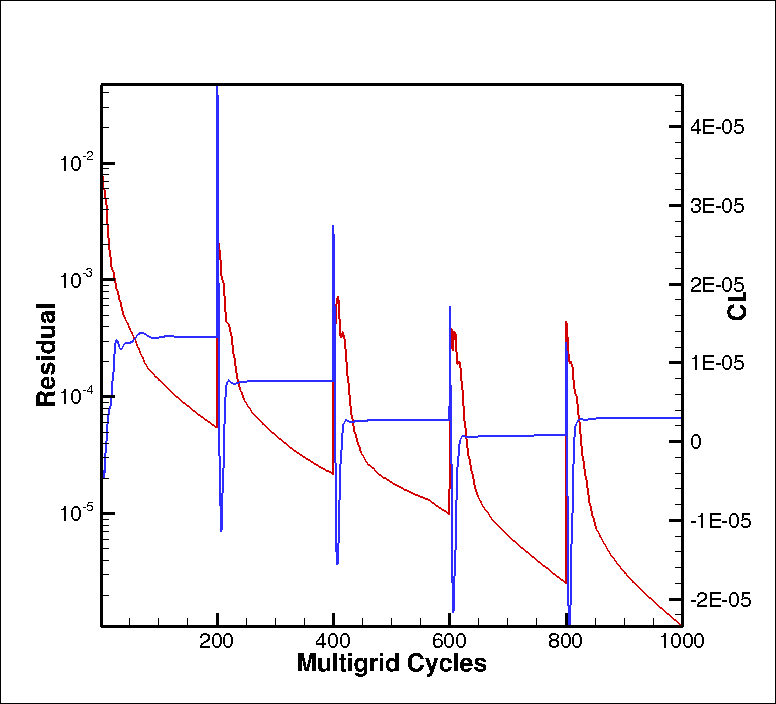} }
  \caption{Illustration of mesh used for rotor test case and nonlinear multigrid convergence history for first 5 time steps}
  \label{fig:rotor1}
\end{figure}

The second test case consists of a time-dependent problem where we examine the efficiency of the various solvers
at converging the non-linear problem arising at each implicit time step using a BDF2 time discretization.
The spatial discretization is identical to that described in the previous case.
The geometry consists of a four-bladed helicopter rotor which is rotated within an oncoming flow to simulate forward flight.
An unstructured mesh of approximately 5 million vertices is used to discretize the flowfield around the rotor,
as illustrated in Figure \ref{fig:rotor1}, using a mixture of tetrahedral elements and highly stretched prismatic elements
in near wall regions, as previously described.
The freestream flow conditions include a Mach number of 0.236 and a Reynolds number
based on the rotor blade chord of 11 million.
To simulate the rotor motion, the entire grid is rotated about the rotor hub axis as a solid body
and a forward shaft tilt angle of 7.31 degrees is prescribed.  These conditions are representative of a high speed forward flight condition
and result in a rotor blade tip Mach number of approximately 0.9, leading to transonic flow on portions of the rotor.
The simulation is initialized using a uniform flowfield corresponding to the freestream conditions
with impulsively started rotor motion using a physical time step corresponding to one degree ($1^o$) of rotation.
At subsequent time steps, the initial condition is taken as the flow solution computed at the end of the previous time step.
This test case is particularly revealing because it highlights the differences between the performance of the nonlinear multigrid
and Newton-Krylov solvers at initial and subsequent time steps in the simulation due to the initial condition.

Figure \ref{fig:rotor1}(c) depicts the convergence rate of the nonlinear multigrid solver for the first 5 time steps.
A fixed number of multigrid cycles are employed at each time step, using a total of 3 multigrid levels.
At each time step, the residual is seen to decrease monotonically and is reduced by approximately 2.5 to 3 orders of magnitude
before being reinitialized for the next time step.  The thrust coefficient, which represents an integrated solution
quantity over the rotor surface (similar to the lift coefficient in the previous case) is also seen to converge reliably to its final value
at each time step over the prescribed number of multigrid cycles. The residual at the start of the first time step
is higher than the initial residual at subsequent time steps, due to the poor initial guess for the impulsively started rotor.
However, the convergence at all time steps is very consistent both in terms of residual decrease and thrust coefficient history.

%----------------------------------------------------------------------
%----------------------------------------------------------------------
%----------------------------------------------------------------------
\begin{figure}[!h]
  \centering
  \subfigure[]{
  \includegraphics[width=0.40\textwidth]{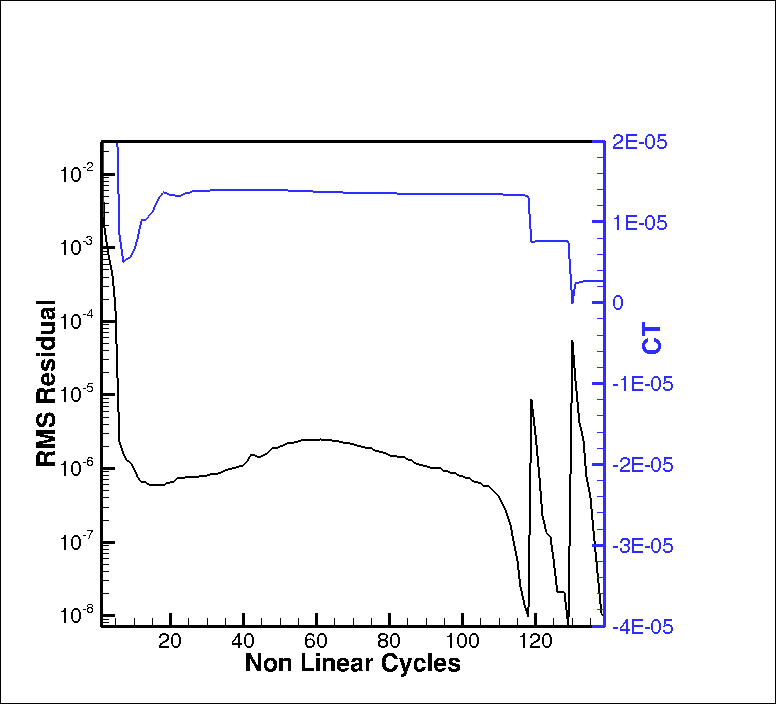}}
  \subfigure[]{
  \includegraphics[width=0.40\textwidth]{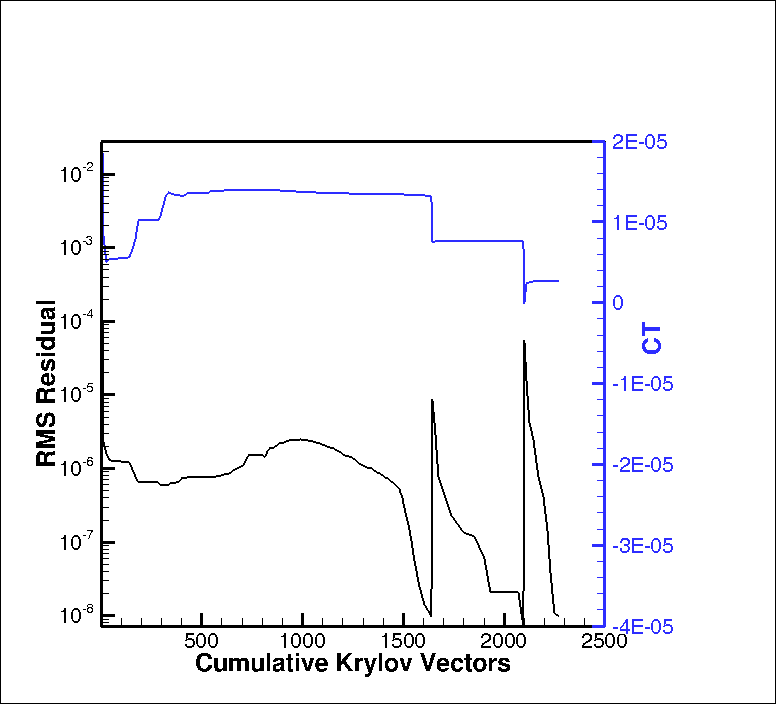}}
  \caption{Convergence for first 3 time steps using unsmoothed Newton-Krylov solver}
  \label{fig:rotor2}
\end{figure}
%----------------------------------------------------------------------
%----------------------------------------------------------------------
\begin{figure}[!h]
  \centering
  \subfigure[]{
  \includegraphics[width=0.40\textwidth]{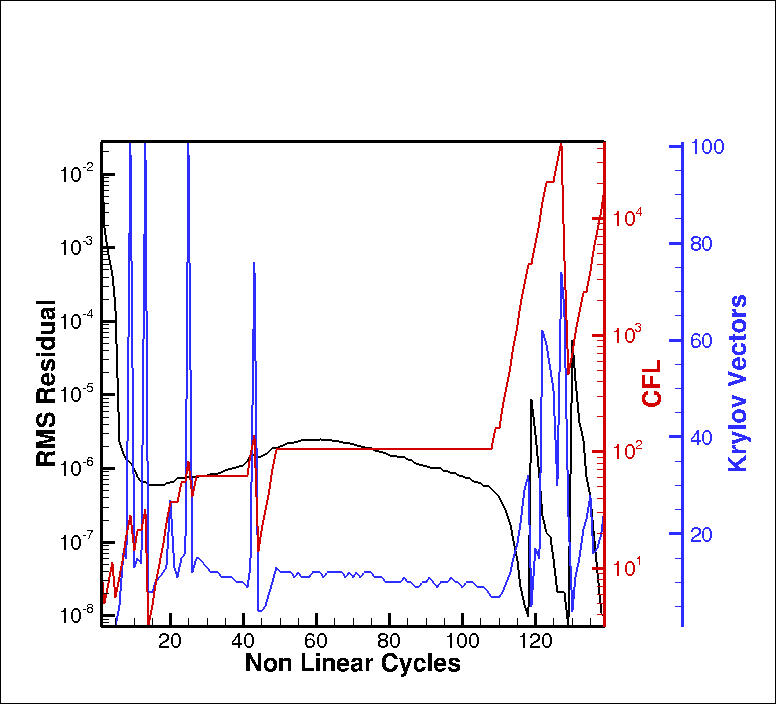}}
  \subfigure[]{
  \includegraphics[width=0.40\textwidth]{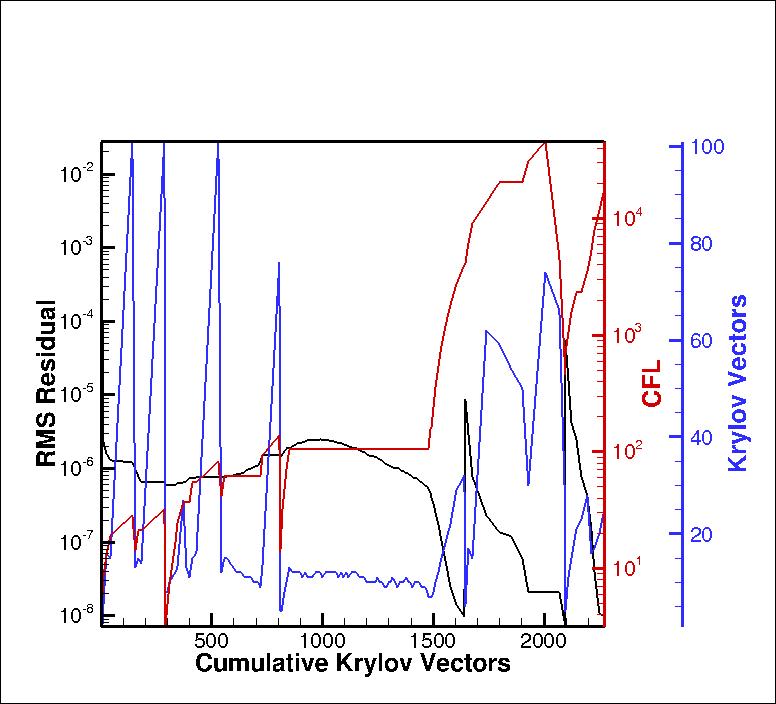}}
  \caption{Convergence history details of pseudo-transient continuation for unsmoothed Newton-Krylov solver}
  \label{fig:rotor3}
\end{figure}
%----------------------------------------------------------------------
%----------------------------------------------------------------------
\begin{figure}[!h]
  \centering
  \subfigure[]{
  \includegraphics[width=0.40\textwidth]{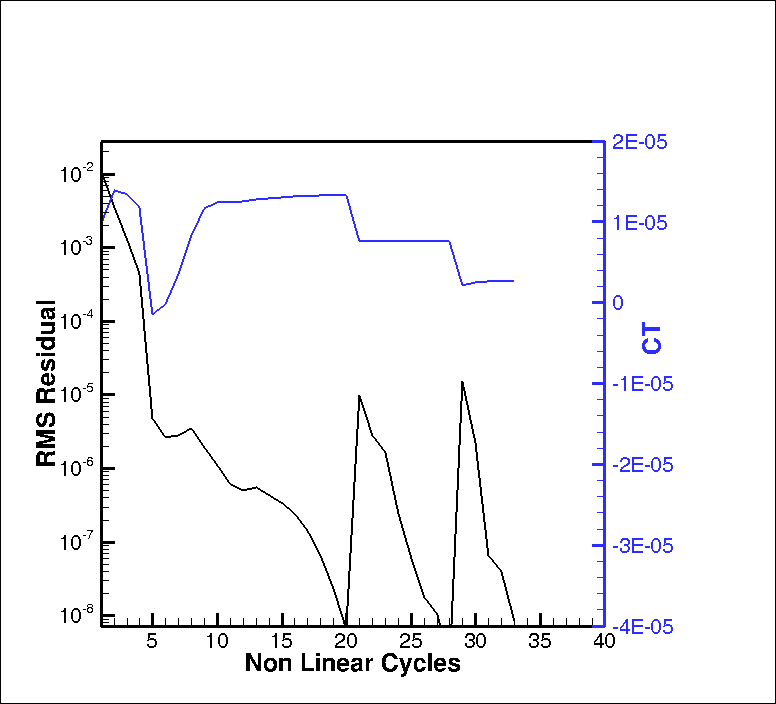}}
  \subfigure[]{
  \includegraphics[width=0.40\textwidth]{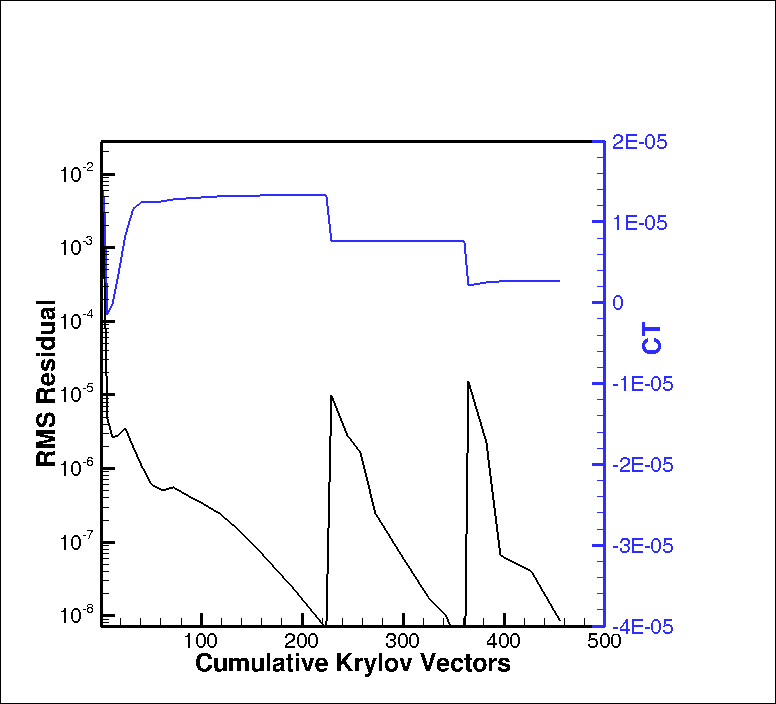}}
  \caption{Convergence for first 3 time steps using smoothed Newton-Krylov solver}
  \label{fig:rotor4}
\end{figure}
%----------------------------------------------------------------------
%----------------------------------------------------------------------
\begin{figure}[!h]
  \centering
  \subfigure[]{
  \includegraphics[width=0.40\textwidth]{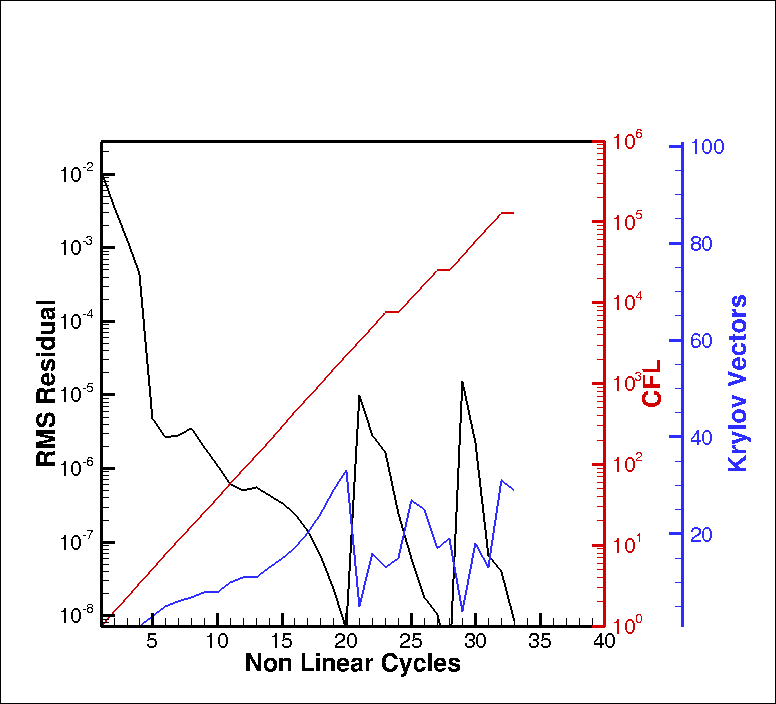}}
  \subfigure[]{
  \includegraphics[width=0.40\textwidth]{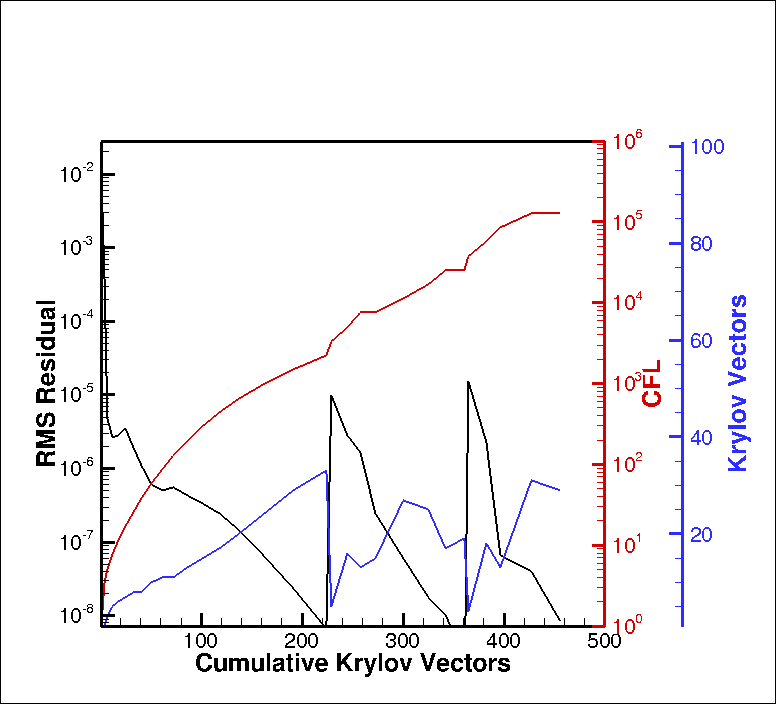}}
  \caption{Convergence history details for first 3 time steps of pseudo-transient continuation for smoothed Newton-Krylov solver}
  \label{fig:rotor5}
\end{figure}
%----------------------------------------------------------------------
%----------------------------------------------------------------------
Figure \ref{fig:rotor2} illustrates the convergence rate of the Newton-Krylov solver for the first 3 time steps
for the same case in terms of nonlinear cycles
and cumulative Krylov vectors. Here the sensitivity of the Newton-Krylov solver to the initial guess at each time step is evident.
On the first time step, this solver requires 118 nonlinear steps and a total of 1639 Krylov vectors
to converge to the prescribed level of 1.e-08, while by the third time step,
the Newton-Krylov solver becomes very effective, rapidly converging the nonlinear system in 10 nonlinear steps and 174 Krylov vectors,
largely due to the fact that the initial guess
provided by the previous time step is close to the domain of quadratic convergence of the Newton method.
Although this behavior is not unexpected, the disparity between the performance of the solver at the first and third time steps
is notable, especially when compared to the performance of the nonlinear multigrid scheme which is similar at all time steps.
Figure \ref{fig:rotor3} provides more details concerning the performance of the pseudo-continuation process in the Newton-Krylov scheme.
For the first time step, growth of the CFL number is impeded by the controller causing a lengthy and costly initial transient
which produces linear systems at various nonlinear cycles that are stiff and require a large number of Krylov vectors,
similarly to the behavior described in the previous test case. However, by the third time step, large CFL values are maintained
and fast convergence is obtained in 10 nonlinear cycles with relatively well behaved linear solver performance using
of the order of 20 Krylov vector per nonlinear step. In essence, the Newton-Krylov scheme is much more costly than the nonlinear multigrid scheme
at converging the initial time step problem, but becomes much more effective than this latter approach for subsequent time steps when
a good initial guess is available.

Figure \ref{fig:rotor4} depicts the convergence history of the residual-smoothed Newton-Krylov approach for this test case.
For this calculation, 5 nonlinear line-preconditioned RK passes were used to construct the smoothing source term,
(similarly to the results shown in Figures \ref{fig:7} and \ref{fig:8} for the previous test case) with all other Newton-Krylov solver settings remaining the same.
From the figure, it is seen that the number of nonlinear cycles required to converge the first time step is reduced from 118 to 20,
with a corresponding drop in the cumulative number of Krylov vectors from 1639 to 224. Furthermore, the disparity between the convergence
rates achieved on the initial and subsequent time steps is much reduced, producing solver behavior which is much more similar
to that observed for the nonlinear multigrid approach.
In this case, the third time step converges in just 5 nonlinear steps requiring a total of 91 Krylov vectors.
Figure \ref{fig:rotor5} provides more details of the continuation process, showing a rapid and monotonic rise
in the CFL number throughout the first time step, and producing linear systems at each nonlinear cycle that are solved
in a relatively small number of Krylov vectors.

\section{Conclusions and Future Work}
In this work he have proposed a residual smoothing approach for addressing the problem of slow initial convergence
in pseudo-transient continuation Newton methods. In contrast with other approaches which,
based on the observed success of local nonlinear solvers, attempt
to break up the problem into smaller more localized nonlinear problems,
we argue that, provided a smooth distribution of residuals can be maintained, the global
nonlinear problem should be able to be advanced efficiently through strong nonlinear transients.
Using a formulation that combines local nonlinear smoothers with a global Newton scheme, significant gains in efficiency
of the overall nonlinear solution process have been demonstrated for realistic CFD problems.
Although the current approach can be interpreted as a strategy for smoothly transitioning from local nonlinear
solvers to a global Newton solver, empirical evidence points to residual smoothing
as the effective mechanism for improving overall nonlinear convergence efficiency.
Further solver efficiency improvements should be possible through the development
of well designed smoothing operators, as well as by disabling the smoothing operations in the final stages
of the continuation process, where the pseudo-transient terms become negligibly small.
Future work will consider using this approach as one component within a broader strategy that seeks
to improve overall Newton-Krylov solver efficiencies through the use of additional continuation techniques
combined with improved linear solver approaches.

\section{Acknowledgments}
This work was partially funded by NASA grant NNX15AU23A under the Tranformational Tools and Technologies ($T^3$) project
and Sandia National Laboratory contract 1852733.
We are grateful for computer time provided by the NCAR-Wyoming Supercomputer Center (NWSC) and by the University of Wyoming Advanced Research Computing Center (ARCC).

\bibliography{/home/dimitri/BIBTEX/UMG}
\bibliographystyle{aiaa}

%------------------------------------------------------------------------------
\end{document}